%% file: main.tex
\documentclass[smallcondensed]{svjour3}
\usepackage[utf8x]{inputenc}
\usepackage{amssymb}
\usepackage{amsmath}
\usepackage[numbers,square]{natbib}
\usepackage{graphics}
\usepackage{graphicx}
\usepackage[lined,boxed,commentsnumbered,linesnumbered]{algorithm2e}
\usepackage{theorem}
\usepackage{tikz}
\usetikzlibrary{calc,3d}
\smartqed

\providecommand{\SetAlgoLined}{\SetLine}
\providecommand{\DontPrintSemicolon}{\dontprintsemicolon}
\newtheorem{Definition}{Definition}

\newtheorem{Theorem}[Definition]{Theorem}

\newtheorem{Lemma}[Definition]{Lemma}
\newtheorem{Corollary}[Definition]{Corollary}

\newtheorem{Proposition}[Definition]{Proposition}
\newenvironment{Proof}{\begin{proof}}{\end{proof}}
\newcommand{\sign}[1]{{\mathrm{sign}(#1)}}

\newcommand{\kerZ}[1]{{\mathrm{ker}(#1)\backslash\{0\}}}
\newcommand{\rg}[1]{{\mathrm{rg}(#1)}}
\newcommand{\homepage}{\texttt{https://www.tu-braunschweig.de/iaa/personal/kruscel}}
\newcommand{\supp}[1]{{\mathrm{supp}(#1)}}

\begin{document}

\title{Computing and Analyzing\\Recoverable Supports for Sparse Reconstruction}
\author{Christian Kruschel \and Dirk A. Lorenz}
\institute{Christian Kruschel \at
  Institute for Analysis and Algebra\\
  TU Braunschweig, 38092 Braunschweig\\
  Tel.: +49-531-3917421\\
  Fax.: +49-531-3917414\\
  \email{c.kruschel@tu-braunschweig.de}
  \and
  Dirk A. Lorenz\at
  Institute for Analysis and Algebra\\
  TU Braunschweig, 38092 Braunschweig\\
  \email{d.lorenz@tu-braunschweig.de}
  }

\date{Received: date / Accepted: date}
 
\maketitle

\begin{abstract}
Designing computational experiments involving $\ell_1$ minimization with linear constraints in a finite-dimensional, real-valued space for receiving a sparse solution with a precise number $k$ of nonzero entries is, in general,  difficult. 
Several conditions were introduced which guarantee that, for small $k$ and for certain matrices, simply placing entries with desired characteristics on a randomly chosen support will produce  vectors which can be recovered by $\ell_1$ minimization.

In this work, we consider the case of large $k$ and propose both a methodology to quickly check whether a given vector is recoverable, and to construct vectors with the largest possible support.
Moreover, we gain new insights in the recoverability in a non-asymptotic regime.
The theoretical results are illustrated with computational experiments.
\end{abstract}

\input{section1_introduction.tex}
\input{section2_existence_and_conditions.tex}
\input{section3_bounds_for_numbers_2.tex}
\input{section4_algorithm.tex}
\input{section5_computational_experiments.tex}

\input{section6_conclusion.tex}

\begin{acknowledgement}
  This work has been funded by the Deutsche Forschungsgemeinschaft
  within the project ``Sparse Exact and Approximate Recovery'' under
  grant LO 1436/3-1.
\end{acknowledgement}

\bibliography{literatur}
\end{document}

%% file: section1_introduction.tex
\section{Introduction}
The difficulty of finding suitable test instances is a serious problem in the field of \textit{Sparse Reconstruction}.
A common and promising method to reconstruct a vector $x^*\in\mathbb{R}^n$ with only a few nonzeros entries from a linear transformation, which is realized by a matrix $A\in\mathbb{R}^{m\times n}$, is performing $\ell_1$ minimization, i.e.
\begin{align}
x^*=\arg\min_y\|y\|_1\mbox{ s.t. } Ay=Ax^*.\label{l1}
\end{align}
This optimization problem was introduced in \cite{ChDoSa89} and is called \textit{Basis Pursuit}. 
Under certain conditions (e.g. see \cite{DoHu01,GrNi03,Tr05}) the vector $x^*$ is also a solution with the smallest number of nonzero entries; a vector $x^*$ with exactly $k$ nonzero entries is called \textit{$k$-sparse}.

A popular method for finding a $k$-sparse vector $x^*$ satisfying \eqref{l1} for a given matrix $A\in\mathbb{R}^{m\times n}$ is to choose an index set $I\subset\{1,\dots,n\}$ with cardinality $k$ and entries $x^*_i, i\in I$, randomly.
For small $k$ this procedure is promising especially if the conditions mentioned above are satisfied, but these conditions require small $k$;
for large $k$ it is more difficult to get suitable $k$-sparse vectors $x^*$.
Besides the question how to compute $x^*$ satisfying \eqref{l1} for a given matrix, we state the question how many different pairs of index sets and signums do exist for a given sparsity $k$.
We aim at partial answers to these questions in a non-asymptotic regime.

We denote by $I=\supp{x^*}$ the support of $x^*$ and its complement by $I^c=\{1,\dots,n\}\backslash I$. With $A_I$ we denote the submatrix of $A$, which columns are indexed by $I$, by $A_I^T$ its transpose, and set $s=\sign{x^*}_I$.
For~(\ref{l1}) to hold, it is necessary and sufficient (cf.~\cite[Theorem 2]{Pl07}) that
\begin{align}
 \exists w\in\mathbb{R}^m : A_I^Tw=s, \|A_{I^c}^Tw\|_\infty<1\mbox{ and }A_I\mbox{ has full rank}.\label{Conditions:l1}
\end{align}
A vector $w$ fulfilling~(\ref{Conditions:l1}) will be called \emph{dual certificate for the support $I$ and sign $s$}.
Condition (\ref{Conditions:l1}) shows that the recoverability of the solution $x^*$ only depends on its support and its signum.
\begin{Definition}
Let $A\in\mathbb{R}^{m\times n}$ and $k\le m\le n$. For $I\subset\{1,\dots,n\}$ and $s\in\{-1,1\}^I$, a pair $(I,s)$ satisfying \eqref{Conditions:l1} is called \textit{Recoverable Support} of $A$. 
If $I$ has the cardinality $k$, a Recoverable Support $(I,s)$ has the \textit{size} $k$.
\end{Definition}

Thus, finding $x^*$ which satisfies \eqref{l1} is equivalent to finding a corresponding Recoverable Support.
For the rest of this paper we will denote the cardinality of a set $I$ with $|I|$ and the $i$-th column of a matrix $A$ with $a_i$.
Moreover we will require $m\le n$ for all ${m\times n}$-matrices.

With a geometrical interpretation of \eqref{Conditions:l1}, new insights to Basis Pursuit, including what kind of matrices can be used and how many Recoverable Supports do exist for a certain size $k$, can be developed.
To that end, consider that $A^Tw$ is a relative interior point of an $(n-|I|)$-dimensional face of the $n$-dimensional hypercube $C^n := [-1,+1]^n$ and assume that the range of $A^T$ is an $m$-dimensional subspace.
Hence, condition \eqref{Conditions:l1} implies the geometrical interpretation that an $m$-dimensional subspace cuts the relative interior of an $(n-|I|)$-dimensional face of $C^n$.
In \cite{Pl07}, the resulting polytope emerging from the intersection of the $m$-dimensional subspace and $C^n$ is considered.
Counting all index sets $I\subset\{1,\dots,k\}$ with $|I|=k$, which satisfy this geometrical interpretation, one can give exact values for the numbers of recoverable vectors for a matrix $A$ and a sparsity $k$.
These values have been estimated in several papers (e.g. \cite{TsDo06,BrDoEl09,DoPeFa10}) through Monte Carlo samplings.
Further this interpretation brings \textit{Sparse Reconstruction} together with the topic \textit{(cross-)sections of a hypercube} in Combinatorial Geometry.

A different geometrical interpretation has been given by Donoho in \cite{Do04} through associating randomly projected $n$-dimensional crosspolytopes with the Basis Pursuit problem, see also the accessible description in~\cite[Section 4.5]{FR13}.
The connection between \textit{Sparse Reconstruction} and the theory of \textit{convex polytopes} gave new insights in both fields. 
Our geometrical interpretation of Recoverable Supports is dual to this approach.
Nonetheless our interpretation delivers additional insights to the questions posed above.

This paper is organzised as follows:
In Section 2 we develop conditions for the existence of Recoverable Supports.
The geometrical aspect around the stated geometrical interpretation will be regarded more carefully in Section 3: a proof for the geometrical interpretation of Recoverable Supports will be given, and exact numbers of Recoverable Supports for certain types of matrices as well as a non-trivial upper bound for these numbers will be stated.
Further we will introduce an algorithm to compute a Recoverable Support of a given matrix and a given size in Section 4.
The theoretical results from these sections will be illustrated by Monte Carlo experiments in Section 5.
Through numercial experiments we additionally provide evidence that checking \eqref{Conditions:l1} is considerably faster than solving Basis Pursuit as a linear program.
In addition, our method stands out from recently done experiments since we can also ensure that a vector is the \emph{unique} solution of Basis Pursuit, without restricting the test problems to a certain class of matrices (e.g. random matrices)
% (e.g. \textit{full spark frames}).
% This is important for applications where a certain vector shall be recovered, as computed tomography.% Ggf. rein mit Referenz auf Sidky, Jörgensen et al.?

%%% Local Variables: 
%%% mode: latex
%%% TeX-master: "main"
%%% End: 

%% file: section2_existence_and_conditions.tex
\section{Existence of and Conditions for Recoverable Supports}
% In the following subsections we establish a partial order on the set of all Recoverable Supports of a given matrix (Section \ref{Se2:Poset}), and give a characteristic that determines if there exists a Recoverable Support of a given matrix (Section \ref{Se2:Ex}).
% This characteristic is so native that it seems only malicious constructed matrices are left behind from the utility of \textit{Sparse Reconstruction}. 
% The necessary and sufficient condition in Section \ref{Se2:Con} is significant for identifying this characteristic.
% Before we establish a partial order and expose maximal elements on the set of all Recoverable Supports of a given matrix in the sense of Zorn's Lemma, labeled as \textit{Maximal Recoverable Supports}.
% The achievements in this section will be used for computing a Recoverable Support in Section 4.

\subsection{Establishing a Partial Order}\label{Se2:Poset}
The condition \eqref{Conditions:l1} for Recoverable Supports rests on two things: The injectivity of the submatrix $A_I$ with $I$ being the support of $x^*$ and the existence of the dual certificate $w\in\mathbb{R}^m$.
The following theorem shows that it is possible to shrink Recoverable Supports and gives conditions when it is possible to obtain a larger Recoverable Support from a given one.
\begin{Theorem}\label{Th:poset}
 Let $A\in\mathbb{R}^{m\times n}$ and let $S_1=(I,s)$ be a Recoverable Support of $A$.
\begin{enumerate}
\item If for $w$ satisfying \eqref{Conditions:l1} there is $y\in\ker{A_I^T}$ satisfying $\|A_{I^c}^T(w+y)\|_\infty=1$ and $A$ restricted to $J = \{i : a_i^T(w+z) = 1\}$ has full rank, then with $t = A_J^T(w+z)$ the pair $S_2=(J,t)$ is a Recoverable Support of $A$ and it holds $I\subset J$.
\item Let $|I|>1$. For any $j_0\in I$ there exists $\tilde{s}\in\{-1,1\}^{I\backslash\{j_0\}}$ with $s_i=\tilde{s}_i$ for all $i\in I\backslash\{j_0\}$, such that the pair $S_3 = (I\backslash\{j_0\},\tilde{s})$ is a Recoverable Support of $A$.
\end{enumerate}
\end{Theorem}
\begin{Proof}
The existence of $y$ for the first statement is obvious and the conclusion that $(J,t)$ is a Recoverable Support follows directly by checking~\eqref{Conditions:l1}. 
For the second statement notice that $A_{I\backslash\{j_0\}}$ has full rank too and secondly that it holds $\ker{A_I^T}\subset\ker{A_{I\backslash\{j_0\}}^T}$.
Hence for $w\in\mathbb{R}^m$ satisfying \eqref{Conditions:l1} there exists $z\in\kerZ{A_{I\backslash\{j_0\}}}$ with $a_{j_0}^Tz \neq 0$.
Choose $\lambda\neq0$ such that $|\lambda|<(1-|a_i^Tw|)/|a_i^Tz|$ for all $i\in I^c$ with $a_i^Tz\neq0$ and
\[\lambda a_{j_0}^Tz\in\left\{\begin{array}{cl}(-2,0)&,\mbox{ if }a_{j_0}^Tw = 1\\(0,2)&,\mbox{ else}\end{array}\right.\]
holds. 
Considering all elements of $A^T(w+\lambda z)$ seperately, it holds
\begin{align*}
|a_{j_0}^Tw + \lambda a_{j_0}^Tz| < 1&,\\
a_i^Tw + \lambda a_i^Tz = a_i^Tw = s_i &\mbox{ for }i\in I\backslash\{j_0\},\\
|a_j^Tw + \lambda a_j^Tz| \le |a_j^Tw| + |\lambda||a_j^Tz| < 1&\mbox{ for }j\in I^c
\end{align*}
by construction. 
Hence with $\tilde{s} = A_{I\backslash\{j_0\}}^T(w+\lambda z)$ the pair $S_2 = (I\backslash\{j_0\},\tilde{s})$ is a Recoverable Support of $A$.
\end{Proof}

The following corollary can be obtained by applying the second statement in Theorem \ref{Th:poset} recursively.
\begin{Corollary}
 Let $A\in\mathbb{R}^{m\times n}$ and $(I,s)$ be a Recoverable Support of $A$. 
Then for any $J\subset I,J\neq\emptyset$, there exists $\tilde{s}\in\mathbb{R}^n$ with $\tilde{s}_J=s_J$, such that the pair $(J,s)$ is a Recoverable Support of $A$.
\end{Corollary}

By using the stated inclusion of Recoverable Supports, a partial order can be obtained through Theorem \ref{Th:poset}:
For Recoverable Supports $S_1 = (I,s), S_2 = (J,\tilde{s})$ with $s_J = \tilde{s}_J$, it is $S_2\le S_1$ if and only if $J\subset I$.
For example, the supports $S_1$, $S_2$ and $S_3$  from Theorem \ref{Th:poset} fulfill $S_3\le S_1\le S_2$.
Moreover, any Recoverable Support can be shrinked and enlarged under the assumption that the respective submatrix is injective.
In other words, the set of all Recoverable Supports form a partially ordered set and may be visualized as a Hasse Diagram.
Further, there exist Recoverable Supports which can not be enlarged, and we call them \emph{Maximal Recoverable Supports}.
Due to the second statement of Theorem \ref{Th:poset}, the Maximal Recoverable Supports determine the full set of all Recoverable Supports.
%In Section 3 we will incorporate Theorem \ref{Th:poset} to establish bounds of the number of Recoverable Supports.

The proof of Theorem \ref{Th:poset} also provides a way to obtain a Recoverable Support if a pair $(I,s)$ satisfies all requirements but having $A_I$ as a full rank matrix. 
\begin{Corollary}
 Let $A\in\mathbb{R}^{m\times n}, I\subset\{1,...,n\}$ and $s\in\{-1,1\}^I$. 
Further let there exist $w\in\mathbb{R}^m$ satisfying $A_I^Tw = s, \|A_{I^c}^Tw\|_\infty<1$. 
If there exists $J\subset I$ such that the submatrix $A_J$ has full rank, then there exists $\tilde{s}\in\{-1,1\}^J$ with $\tilde{s}_J = s_J$ such that $(J,\tilde{s})$ is a Recoverable Support of $A$.
\end{Corollary}

\subsection{Sufficient and Necessary Condition}\label{Se2:Con}
Similar to Section \ref{Se2:Poset}, we will consider dual certificates to establish a sufficient and necessary condition for a pair $(I,s)$ being a Recoverable Support of a given matrix.
For this purpose, we introduce the pseudo-inverse $(A_I^T)^\dagger$ of $A_I^T$.
The following theorem und its corollary are an extension of Fuchs' sufficient condition in \cite{Fu04}.
\begin{Theorem}\label{Th:SuffNessCond}
Let $A\in\mathbb{R}^{m\times n}, I\subset\{1,...,n\}$ and $s\in\{-1,1\}^I$. 
Then $(I,s)$ is a Recoverable Support of $A$ if and only if $A_I$ has full rank and there exists $y\in\ker{(A_I^T)}$ such that
\begin{align*}
\|A_{I^c}^T(A_I^{T})^\dagger s+ A_{I^c}^Ty\|_\infty <1.
\end{align*}
\end{Theorem}
\begin{Proof}
If $(I,s)$ is a Recoverable Support, then $A_I$ has full rank and there exists $w\in\mathbb{R}^m$ such that $A_I^Tw=s$. 
For $\tilde{y}\in\mbox{ker}(A_I^T)$ the solution has the general representation $w = (A_I^{T})^\dagger s+\tilde{y}$.
Since there exists at least one $w$ satisfying $\|A_{I^c}^Tw\|_\infty<1$, there exists $y\in\mbox{ker}(A_I^T)$ proving the stated inequality.

Further for $y\in\mbox{ker}(A_I^T)$ consider $w= (A_I^{T})^\dagger s + y$. 
Since $A_I$ has full rank, $A_I^T$ has linearly independent rows, so $A_I^Tw = s$ holds as well as $\|A_{I^c}^Tw\|_\infty<1$.
\end{Proof}

Note that a conclusion of Theorem \ref{Th:SuffNessCond} is that
\[A_I\mbox{ has full rank and }\|A_{I^c}^T(A_I^{T})^\dagger s\|_\infty <1\]
is a sufficient condition for $(I,s)$ being a Recoverable Support of $A$ by choosing $y=0$. 
For full rank matrices with $|I|=\mbox{rank}(A)$ this is also a necessary condition using the inverse $A_I^{-T}$ of $A_I^T$.
\begin{Corollary}\label{Co:SuffNessCondFullRank}
Let $A\in\mathbb{R}^{m\times n}$ have a full rank, $I\subset\{1,...,n\}$ with $|I| = m$ and $s\in\{-1,1\}^I$. 
Then $A_I$ is invertible and it holds $\|A_{I^c}^TA_I^{-T}s\|_\infty <1$ if and only if $(I,s)$ is a Maximal Recoverable Support of $A$.
\end{Corollary}

% Theorem \ref{Th:SuffNessCond} will be applied in the coming subsection.
% Necessary to that end is that the pseudo-inverse of a vector $x$ is interpreted as $\|x\|^{-2}x^T$.

%\subsection{Existence of Recoverable Supports}\label{Se2:Ex}
We close this section with the observation that there exist matrices which do not possess any Recoverable Support.
The following theorem characterizes these matrices. %Note that the pseudo-inverse of a vector $x$ is $\|x\|^{-2}x^T$.
\begin{Theorem}\label{Th:Existence}
Let $A\in\mathbb{R}^{m\times n}$ and $k\in\{1,...,n\}$ such that for all $j\neq k$ holds $\|a_j\|\le\|a_k\|$. 
Then for $s\in\{-1,+1\}$ the pair $(\{k\},s)$ is a Recoverable Support of $A$ if and only if for any $j\neq k$ with $\|a_j\|=\|a_k\|$ it holds $a_j\neq a_k$.
\end{Theorem}

\begin{Proof}
Let $(\{k\},s)$ be a Recoverable Support of $A$ and without loss of generality let $s=+1$. 
Assuming for $j\neq k$ it holds $a_k=a_j$, then for all $y\bot a_k$ it holds
\[\left|\|a_k\|^{-2}a_j^Ta_k+a_j^Ty\right|=1\]
which is a contradiction to Theorem \ref{Th:SuffNessCond}.

For the converse implication let $a_j\neq a_k$ with $\|a_j\| = \|a_k\|$. 
With $w = \|a_k\|^{-2}a_k$ it holds 
\[|a_k^Tw| = 1\mbox{ and }|a_j^Tw| = \frac{|a_j^Ta_k|}{\|a_k\|^2} < \frac{\|a_j\|}{\|a_k\|} = 1\]
by applying Cauchy-Schwartz inequality.
Further for any $a_i$ satisfying $\|a_i\|<\|a_k\|$ the inequality $|a_i^Tw|<1$ holds.
Trivially $A_{\{k\}}$ has full rank, so with $s = a_k^Tw$ it holds that the pair $(\{k\},s)$ is a Recoverable Support of $A$.
\end{Proof}

Hence, every matrix for which the largest column does not appear multiple times possesses a Recoverable Support.
Moreover, Theorem \ref{Th:Existence} will be useful as a starting point for the algorithm in Section 4.
%%% Local Variables: 
%%% mode: latex
%%% TeX-master: "main"
%%% End: 

%% file: section3_bounds_for_numbers_2.tex
\section{Geometrical Interpretation and Number of Recoverable Supports}
In this section, we deal with the geometrical interpretation of Recoverable Supports presented in Section 1 and its implications on their number.
In the end of this section, we further derive a non-trivial, but heuristic upper bound on this number, which is, as far as we know, new.

\begin{Definition}\label{Def:Bounds}
 For $A\in\mathbb{R}^{m\times n}$ the number $\Lambda(A,k)$ is defined as
\[\Lambda(A,k) := |\{(I,s) : (I,s)\mbox{ is a Recoverable Support of }A\mbox{ with size }k\}|.\]
Further let $\Xi(m,n,k)$ be defined as the maximum of $\Lambda$ over all matrices, i.e.
\[\Xi(m,n,k) := \max\{\Lambda(A,k) : A\in\mathbb{R}^{m\times n}\}.\]
\end{Definition}%
For some triples $(m,n,k)$, the values for $\Lambda$ and $\Xi$ will be derived in Sections \ref{Se3:Do} and \ref{Se:ValLam}.
Prior, we briefly sketch some basics on \textit{convex polytopes} in the next section.
% Section \ref{Se3:Pr} and give the mentioned geometrical interpreation of Recoverable Supports in Section \ref{Se3:GeoInt}.
% The approached geometrical interpretation established in \cite{Do04} and its relationship to Theorem \ref{Th:GeometricalInterpretL1} will be discussed in Section \ref{Se3:Do}.

\subsection{Preliminaries}\label{Se3:Pr}
Let $x_1,...,x_m\in\mathbb{R}^n$, then its convex hull $P=\mbox{conv}(x_1,...,x_n)$ is called a \textit{polytope}. 
The \textit{dimension} of a polytope is the dimension of its affine hull; a polytope with dimension $d$ is called \textit{$d$-polytope}. 
We call $P$ \textit{centrally-symmetric} if for all $x\in P$ it holds $-x\in P$.
For $\lambda\in\mathbb{R}^n$ and $c\in\mathbb{R}$ we define the hyperplane $H_{\lambda,c}=\{x:\lambda^Tx=c\}$.
Further the intersection $F=H_{\lambda,c}\cap P$ is called a \textit{face} of $P$ if $\lambda^Tx < c$ holds for all $x\notin F$.
A face of $P$ is also a polytope; more general, any intersection of a polytope with an affine subspace is a polytope. 
The set of all $k$-dimensional faces of $P$ is denoted as $\mathcal{F}_k(P)$. 
A centrally symmetric polytope is called $k$-neigborly if any set of $k+1$ vertices of $P$, not including an antipodal pair, spans a face of $P$.

A face $F$ of the hypercube $C^n := [-1,+1]^n$ is uniquely determined by a pair $(I,s)$ consisting of an index set $I\subset\{1,...,n\}$ and $s\in\{-1,+1\}^I$:
With $(I,s)$ choose $\lambda\in\mathbb{R}^n$ through $\lambda_I=s, \lambda_j = 0$ if $j\notin I$.
We see that for any $y\in\mathbb{R}^n$ with $\lambda^Ty > n-|I|$ it holds $y\notin C^n$.
Hence, it holds that $F=H_{\lambda,(n-|I|)}\cap C^n$ is an $|I|$-dimensional face of $C^n$.
For $F\subset C^n$ we note the following equivalence:
\begin{align}
\begin{array}{c}
 F\in\mathcal{F}_k(C^n)\\
\Leftrightarrow\\
\exists!I\subset\{1,...,n\},|I|=n-k\mbox{ }\forall v,w\in F: v_I\in\{-1,1\}^I, v_I = w_I.
\end{array}
 \label{Eq:FaceIndex}
\end{align}
With $I(F)$ we denote the unique subset of $\{1,...,n\}$ determined by $F\in\mathcal{F}_k(C^n)$.
Since the equivalence also holds for subsets $V\subset F$, we also use $I(V)$ to denote the unique subset.
We collect these observations in the next lemma.
\begin{Lemma}\label{Lemma:faceofcube}
 For $F\in\mathcal{F}_k(C^n)$ there exists $\lambda\in\mathbb{R}^n$ such that
\[F = \{x\in C^n:\lambda^Tx = n-k\}\mbox{ and }C^n\backslash F = \{y\in C^n : \lambda^Ty<n-k\}.\]
\end{Lemma}
%\begin{Proof}
%Consider for $I = I(F)$ the vector $\lambda\in\mathbb{R}^n$ with $\lambda_I = v_I$ for all $v\in F$ and $\lambda_j = 0$ if $j\in I^c$. 
%It holds $\lambda\in F$.
%
%If $x\in C^n$ satisfies $\lambda^Tx = n-k$ then $n-k = \lambda^Tx = \sum_{i\in I}\pm1_i x_i$ and since $|I|=n-k$ it holds $x_i = \pm1_i$ and further $x_I = \lambda_I$. On the other hand, if $x\in F$, then $\lambda_I=v_I$ and $\lambda^Tx = |I|$.
%
%If $y\in C^n\backslash F$ and $\lambda^Ty > n-k$ then there exists $i\in I$ with $|y_i|>1$. Further for $y\in C^n$ with $\lambda^Ty<n-k$ there exists $i\in I$ with $|y_i|<1$ which implies $y_I\neq\lambda_I$ and $y\in C^n\backslash F$.
%\end{Proof}
On the basis of Lemma \ref{Lemma:faceofcube}, we identify the relative interior of a face $F$ with $\mbox{relint}(F)=\{x\in F : |x_i|<1, i\notin I(F)\}$.

For an extensive overview in the field of convex polytopes, we refer to the books by Gr\"unbaum \cite{Gr03} or Ziegler \cite{Zi95}.

Finally, we have all tools for proving the geometrical interpreation of Recoverable Supports suggested in Section 1.

\subsection{Geometrical Interpretation of Recoverable Supports}\label{Se3:GeoInt}
With the introduced notation we will prove the following theorem. Note that the results are similar to the interpretation in \cite{Pl07}.
\begin{Theorem}\label{Th:GeometricalInterpretL1}
Let $A\in\mathbb{R}^{m\times n}$ have rank $l$ and let $k\le l$. 
Then the following statements are equivalent:
\begin{enumerate}
\item\label{EquiTh1} There exists a Recoverable Support of $A$ with size $k$.
\item\label{EquiTh2} There exists $F\in\mathcal{F}_{n-k}(C^n)$ such that $\mbox{relint}(F)\cap\rg{A^T}\neq\emptyset$ and $A_{I(F)}^T$ has full rank.
\item\label{EquiTh3} There exists $V\in\mathcal{F}_{l-k}(C^n\cap\rg{A^T})$ and $v\in V$ with $\|v_{I^c}\|_\infty<1$ and $A_{I}$ has full rank for $I:=I(V)$.
\end{enumerate}
\end{Theorem}
\begin{Proof}
First we state for any subset $I\subset\{1,...,n\}$ with $|I|\le l$ that $A_I$ has full rank if and only if $A_I^T$ has full rank.

$\eqref{EquiTh1}\Rightarrow\eqref{EquiTh2}:$ 
Let $(I,s)$ be a Recoverable Support of $A$ with size $k$.
Choose $\lambda\in\mathbb{R}^n$ with $\lambda_I = s$ and $\lambda_j = 0$ for $j\in I^c$ and consider $F := \{x\in C^n : \lambda^Tx = k\}$. 
It holds $F\in\mathcal{F}_{n-k}(C^n)$.
By assumption there is $W\in\mathbb{R}^m$ such that $A^Tw\in\mbox{relint}(F)$, hence  $\mbox{relint}(F)\cap\rg{A^T}\neq\emptyset$.

$\eqref{EquiTh2}\Rightarrow\eqref{EquiTh3}:$ 
Let $I=I(F)$ and choose $\lambda\in F$ with $\lambda_ = 0$ for $j\notin I$. Then $V=\{y\in P : \lambda^Ty=k\}$ is a face of $P=C^n\cap\rg{A^T}$ and further $V\subset F$.
Hence $I=I(V)$ and there is $v\in V$ with $\|v_{I^c}\|_\infty<1$.
Since $P$ is an $l$-polytope, it holds $V\in\mathcal{F}_{l-k}(P)$.

$\eqref{EquiTh3}\Rightarrow\eqref{EquiTh1}:$ 
Let $I=I(V)$.
There is $w\in\mathbb{R}^m$ such that $A^Tw = v$ and further $\|A_{I^c}^Tw\|_\infty<1$.
Hence, the pair $(I,v_I)$ is a Recoverable Support of $A$ with size $k$.
\end{Proof}

Theorem \ref{Th:GeometricalInterpretL1} partitions solutions of \eqref{l1} into equivalence classes separated into faces of $C^n$ with different dimensions.
For the rest of this section, we will use the notation of each polytope used in Theorem \ref{Th:GeometricalInterpretL1} and $P := C^n\cap\rg{A^T}$.
A first consequence of the latter theorem gives an equivalent expression of Definition \ref{Def:Bounds}: For $A\in\mathbb{R}^{m\times n}$ with rank $l$ and $k\le l$ it is $\Lambda(A,k) = |\mathcal{F}_{l-k}(P)|$.

Further the second statement from Theorem \ref{Th:poset} delivers the following corollary.
\begin{Corollary}\label{Co:GeometricalInterpretL1}
Let $A\in\mathbb{R}^{m\times n}$ have rank $l$. Then the polytope $P = C^n\cap\rg{A^T}$ is $l$-dimensional, centrally-symmetric, and simple, i.e. any vertex of $P$ is adjacenced by $l$ edges.
\end{Corollary}

With Corollary \ref{Co:GeometricalInterpretL1} we can link \textit{Sparse Reconstruction} to simple, centrally-symmetric polytopes. 
Further with the two representations of the geometrical interpretation given by Theorem \ref{Th:GeometricalInterpretL1} we can involve the results from the field \textit{(cross-)sections of a hypercube} from Combinatorial Geometry.
This will be done in Section \ref{Se:ValLam} and \ref{Se3:Bounds}.

\subsection{Geometrical Interpretation of Basis Pursuit by Donoho}\label{Se3:Do}
In this subsection we briefly present the geometrical interpretation of Basis Pursuit by Donoho \cite{Do04,Do04b}.

With the crosspolytope $\mathcal{C} = \{x\in\mathbb{R}^n : \|x\|_1\le1\}$ and the projection operator $A$, we consider the projected crosspolytope $A\mathcal{C} = \{Ax : x\in\mathcal{C}\}$ and further the following theorem.
\begin{Theorem}\cite[Theorem 1]{Do04}\label{Th:GeometricalInterpretL1Do}
Let $A\in\mathbb{R}^{m\times n}$. These two statements are equivalent:
\begin{itemize}
 \item The polytope $A\mathcal{C}$ has $2n$ vertices and is $k$-neighborly.
 \item Any $k$-sparse vector solves Basis Pursuit uniquely.
\end{itemize}
\end{Theorem}
Theorem \ref{Th:GeometricalInterpretL1Do} connects Sparse Reconstrunction with projected crosspolytopes.
Thus, one can apply results from convex polytopes like the following necessary condition taken from~\cite{Do04} which is based on \cite{MuSh68}.
\begin{Corollary}\label{Co:NecK}
 Let $A\in\mathbb{R}^{m\times n}$ with $2 < m \le n-2$. If any $k$-sparse vector $x^*$ solves \eqref{l1} then $k\le\lfloor(m+1)/3\rfloor$.
\end{Corollary}

Further in \cite{Do04b} tools from \cite{AfSc92} are used to count the faces of randomly-projected crosspolytopes.
Considering that any preimage of a face of $A\mathcal{C}$ is a face of $\mathcal{C}$ (e.g. see \cite[Theorem 7.10]{Zi95}), the following lemma connects the property of $k$-neigborliness of $A\mathcal{C}$ and $\mathcal{C}$.
We need the term \textit{$k$-simplex} describing a polytope with $k+1$ vertices.
\begin{Lemma}\cite[Lemma 2.1]{Do04b}\label{Le:FacesKNeighbor}
Let $A$ be a projection and $P = A\mathcal{C}$ such that for $k\in\mathbb{N}$ it holds $|\mathcal{F}_i(P)| = |\mathcal{F}_i(\mathcal{C})|$ for $i=1,...,k-1$.
Then any $F\in\mathcal{F}_l(P)$ is an $l$-simplex for $l=0,...,k-1$ and $P$ is $k$-neighborly.
\end{Lemma}

Hence, with Lemma \ref{Le:FacesKNeighbor} one can say rakishly that if we are losing faces through the projection, Basis Pursuit loses the power of reconstructing sparse vectors.
Moreover, there exist explicit functions $\rho_N,\rho_F:(0,1]\rightarrow[0,1]$ (cf. \cite[Section 3]{Do04b}) such that the following theorems hold.
\begin{Theorem}\cite[Theorem 1]{Do04b}\label{Th:PropNeighbStrong}
 Let $\rho<\rho_N(\delta)$ and $A:\mathbb{R}^n\rightarrow\mathbb{R}^m$ a uniformly-distributed random projection with $m\ge\delta n$. 
Then
\[Prob\{|\mathcal{F}_l(\mathcal{C})|=|\mathcal{F}_l(A\mathcal{C})|, l=0,...,\lfloor\rho m\rfloor\}\rightarrow1\mbox{ as }n\rightarrow\infty.\]
\end{Theorem}
% A weaker notion of neighborliness ic contained in Theorem \ref{Th:PropNeighbWeak}.
% There, an approximate equality of face numbers in the proportional-dimensional case is studied.
\begin{Theorem}\cite[Theorem 2]{Do04b}\label{Th:PropNeighbWeak}
Let $m\sim\delta n$ and $A:\mathbb{R}^n\rightarrow\mathbb{R}^m$ be a uniform random projection. 
Then for $k$ with $k/m\sim\rho, \rho<\rho_F(\delta)$ it holds
\[|\mathcal{F}_k(A\mathcal{C})| = |\mathcal{F}_k(\mathcal{C})|(1+o(1)).\] 
\end{Theorem}
\begin{figure}
\begin{center}\includegraphics[width=80mm]{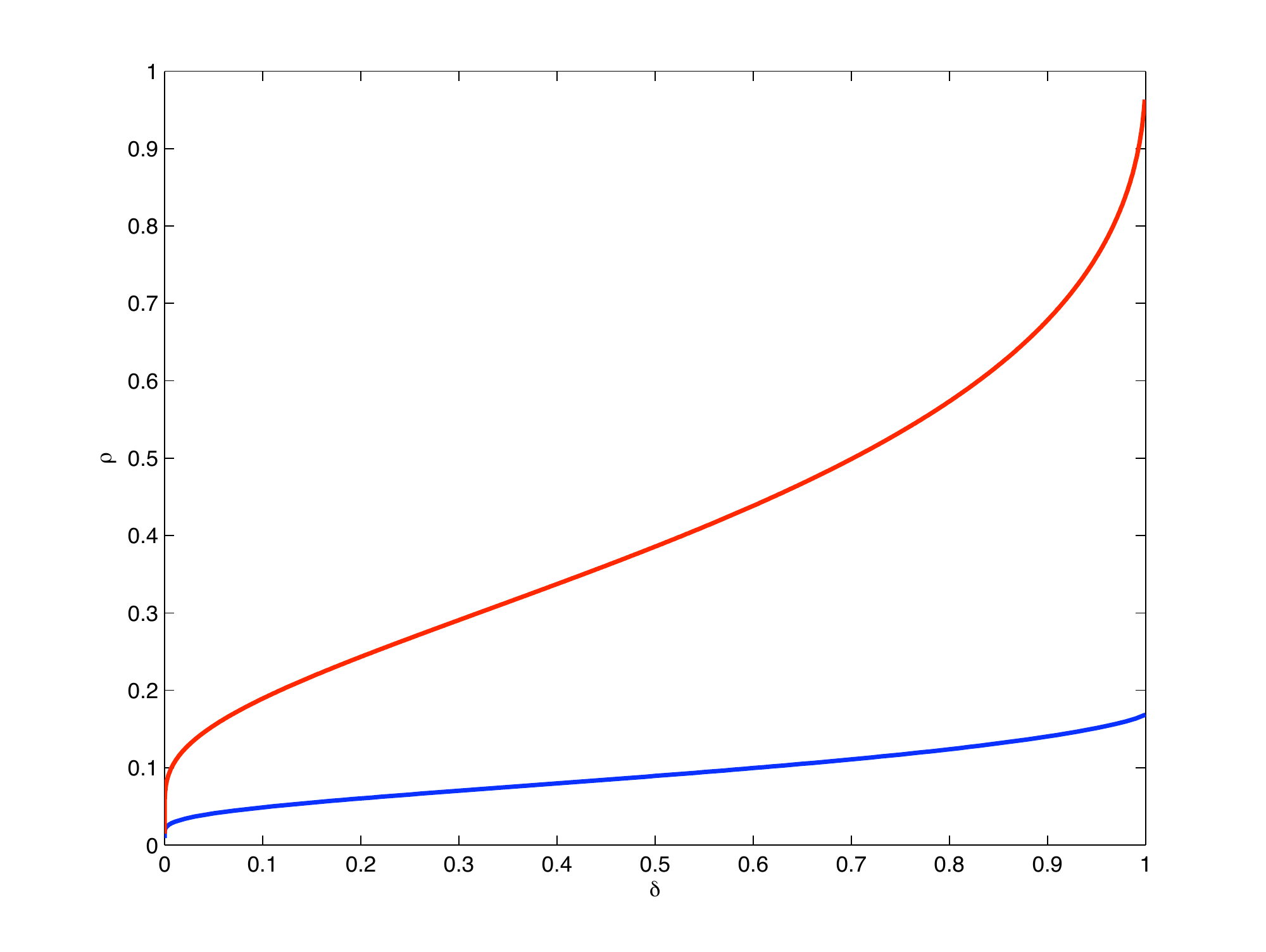}\end{center}
\label{Fig:PhaseTransition}
\caption{Functions $\rho_N$ (blue) and $\rho_F$ (red) in Theorems \ref{Th:PropNeighbStrong} and \ref{Th:PropNeighbWeak}.}
\end{figure}
The functions $\rho_N,\rho_F$ are displayed in Figure \ref{Fig:PhaseTransition} and are known in the context of \textit{Phase Transitions}~\cite{DoTa09}.
Theorem \ref{Th:PropNeighbStrong} implies that for large $m$ and $n$ tending to infinity, with high probability any $\lfloor\rho m\rfloor$-sparse vector $x^*$ is recoverable.
Donoho states \cite[Section 1.5]{Do04b} that the result in Theorem \ref{Th:PropNeighbWeak} can be seen ``as a weak kind of neighborliness [...] in which the overwhelming majority of (rather than all) $k$-tuples span $(k-1)$-faces''.
Further he remarks that this result is ``sharp in the sense that for sequences with [$k/m\sim\rho>\rho_F(\delta)$], we do not have the approximate equality''.
An additional result \cite[Theorem 4]{Do04b} is the limit value consideration
\[\lim\limits_{\delta\rightarrow1}\rho_F(\delta) = 1.\]
This value combined with Theorem \ref{Th:PropNeighbWeak} implies that for $\delta\rightarrow1$ and $n\rightarrow\infty$ \textit{almost} all vectors $x^*$ can be recovered through \eqref{l1} since the number $|\mathcal{F}_k(A\mathcal{C})|$ tends to concentrate near its upper bound value $2^{k+1}\binom{n}{k+1}$.

Taking up our geometrical perspective, we introduce the \textit{polar set} $K^*$ of $K\subset\mathbb{R}^m$ as
\[K^* := \{w\in\mathbb{R}^m : x^Tw\le1\mbox{ for all }x\in K\}\]
and see with 
\[(A\mathcal{C})^* = \{w\in\mathbb{R}^m : |a_i^Tw|\le 1, i=1,...,n\} = \rg{A^T}\cap C^n\]
that the projected cross-polytope $A\mathcal{C}$ is dual to $P$ in Theorem \ref{Th:GeometricalInterpretL1}, see also \cite{Pl07}.
Hence, our approach simply differs that we additionally consider unique solutions of Basis Pursuit.

For further considerations, we denote the cross-section of an $m$-dimensional subspace $K$ of $\mathbb{R}^n$ and $C^n$ as \textit{regular} if $K$ has no point in common with any $(n-m-1)$-dimensional face of $C^n$.
The second statement in Theorem \ref{Th:GeometricalInterpretL1} connects regular cross-sections of the hypercube to Recoverable Supports.
In general, we can not assume that the sections occuring through regarding the range of $A^T$ are regular but we still can use some basic result from literature and connect them to \textit{Sparse Reconstruction}.
This is done in Section \ref{Se:ValLam} and \ref{Se3:Bounds}.

\subsection{Values for $\Lambda$}\label{Se:ValLam}
In this subsection, we give some values of $\Lambda$ for specified matrices and sizes of their Recoverable Supports. 
In general, the polytope $P=\rg{A^T}\cap C^n$ is not a regular cross-section.
Thus, the already difficult problem of counting $k$-faces of a (simple) polytope becomes even more difficult counting only all $k$-faces of $P$ intersecting with $n-m+k$-faces of $C^n$ in case of full rank matrices. 
Different from $\Xi$ (cf. Section \ref{Se3:Bounds}), using past results for a lower bound of $\Lambda$ over all $m\times n$-matrices is, as far as we know, only possible under certain assumptions, as the following corollary states.
\begin{Corollary}\label{Co:LowerBound}
 Let $A\in\mathbb{R}^{m\times n}$ with rank $l$ and assume $\rg{A^T}\cap C^n$ is a regular cross-section.
Then
\[\Lambda(A,l)\ge2^l.\]
\end{Corollary}
\begin{Proof}
 The result follows from Statement 3 of Theorem \ref{Th:GeometricalInterpretL1} and \cite[Corollary 2]{BaLo82}.
\end{Proof}

With the same assumptions, Euler's relation \cite{Po93,Po99} and Steinitz' characterization for $3$-polytopes \cite{St06} can be applied, but the practicability is limited since for every matrix the regularity of its corresponding cross-section has to be checked.
Considering the cross-section as a simple polytope delivers a different lower bound, which is only dependent on the value $\Lambda(A,1)$.
\begin{Corollary}\label{Co:Co:LowerBound_Simple}
 Let $A\in\mathbb{R}^{m\times n}$ with rank $l$. Then
\[\Lambda(A,l)\ge(l-1)\Lambda(A,1) - (l+1)(l-2).\]
\end{Corollary}
\begin{Proof}
 Combining \cite[Theorem 1]{Ba71} and Corollary \ref{Co:GeometricalInterpretL1} proves the result.
\end{Proof}

Note that Corollary \ref{Co:Co:LowerBound_Simple} provides a lower bound on the number of Recoverable Supports of a matrix if the number of Recoverable Supports of size one is known. However, there are no more than $2n$ possibilities and these can be checked easily for any matrix.

For the rest of this section we consider two types of matrices: \textit{Equiangular tight frames} and \textit{Gaussian matrices}.
The term \textit{equiangular tight frame} will be dwelled on later; a Gaussian matrix means that its entries are independant and standard normally distributed random variables, i.e. having mean zero and variance one.

First we consider Gaussian matrices and regard the work of Lonke in \cite{Lo00}.
With $\mbox{erf}$ we denote the \textit{Gauss Error function} and $\mathcal{E}(Z)$ describes the expected value of $Z$. 
\begin{Corollary}\label{Co:EstRandVertices}
Let $A\in\mathbb{R}^{m\times n}$ be a randomly drawn Gaussian matrix.
Then
\[\mathcal{E}(\Lambda(A,m)) = 2^m\binom{n}{m}\sqrt{\frac{2m}\pi}\int_0^\infty e^{-mt^2/2}\left[\mathrm{erf}\left(\frac t{\sqrt{2}}\right)\right]^{n-m}dt.\]
Further it holds
\begin{align}
\mathcal{E}(\Lambda(A,m))\ge\binom{n}{n-m}2^n\left(\frac 1\pi\arctan\frac1{\sqrt{m}}\right)^{n-m}\label{Eq:EstRandVertices}
\end{align}
where for $m = n-1$ equality holds.
\end{Corollary}

\begin{Proof}
 The result follows from \cite[Proposition 2.2, Proposition 2.5]{Lo00} and the second statement of Theorem \ref{Th:GeometricalInterpretL1}.
\end{Proof}

In Section 5 we will match \eqref{Eq:EstRandVertices} with Monte-Carlo samplings.
Additionally, Lonke delivers an asympotic behavior for sizes $k\neq m$.
\begin{Corollary}
Let $A\in\mathbb{R}^{m\times n}$ be a randomly drawn Gaussian matrix. 
Then for $k\neq m$ it holds
\[\lim\limits_{n\rightarrow\infty}\mathcal{E}(\Lambda(A,k))(2n)^{-k}k! = 1.\]
\end{Corollary}

\begin{Proof}
 Combining \cite[Corollary 3.4]{Lo00} and the second statement of Theorem \ref{Th:GeometricalInterpretL1} proves the assertion.
\end{Proof}

As Lonke says \cite[Section 3]{Lo00}, the value $\Lambda(A,k)$ ``tends to concentrate near the value [$2^{k}\binom{n}{k}$], which bounds it from above'' (cf. the statement of Donoho \cite{Do04b} as an implication of Theorem \ref{Th:PropNeighbWeak}).

\begin{figure}
\begin{center}\includegraphics[width=100mm]{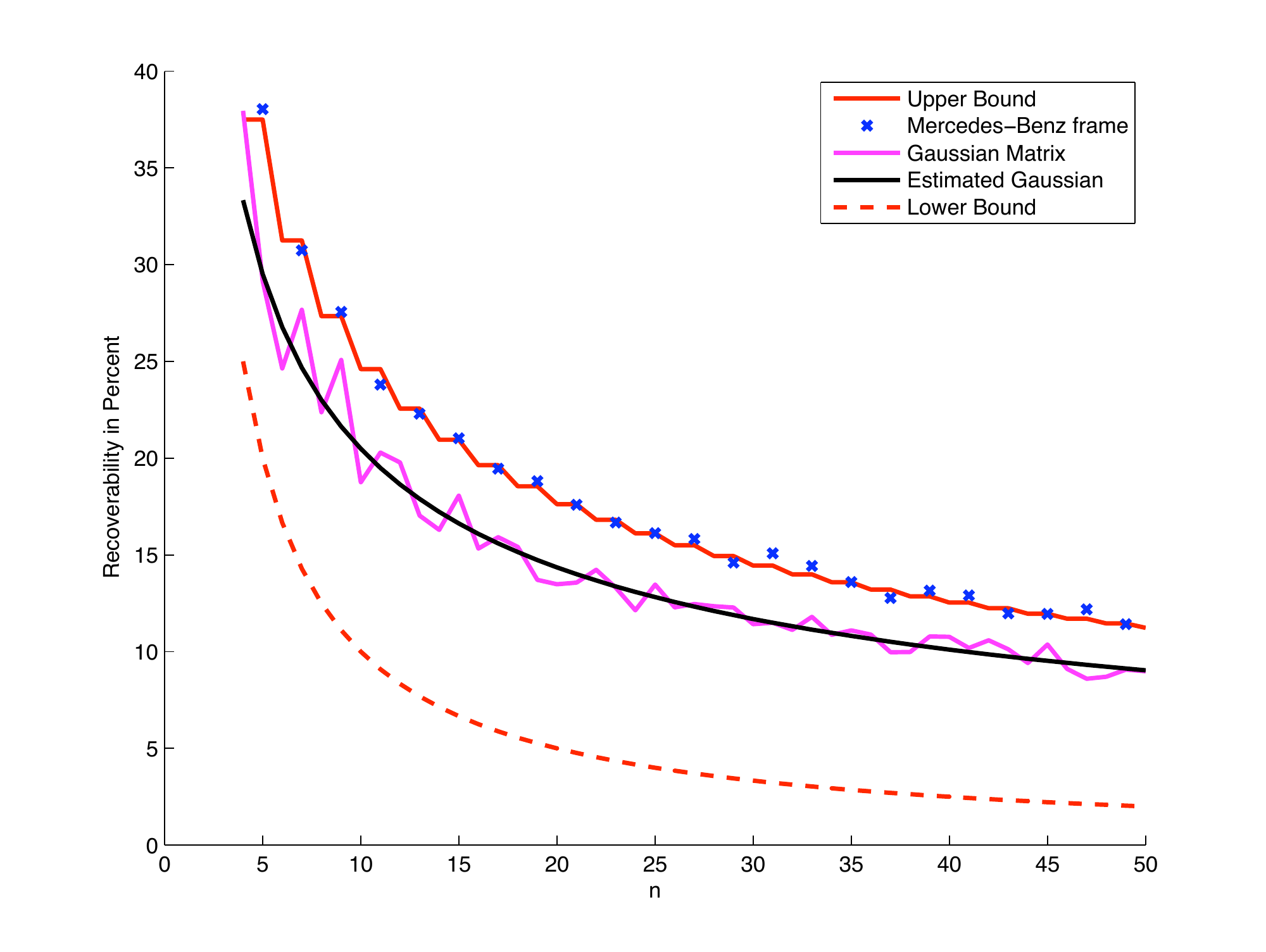}\end{center}
\label{Fig:TheoResults}
\caption{Monte Carlo Sampling versus Theorectical Results in Section \ref{Se:ValLam} and \ref{Se3:Bounds}.
Results from Monte Carlo experiments for Mercedes-Benz frame (blue) and Gaussian matrix (magenta) of the size $(n-1)\times n)$.
For any $n\ge4$, one thousand pairs $(I,s)$ with $I\subset\{1,...,n\}, |I|=n-1, s\in\{-1,+1\}^I$ where taken randomly and tested whether $(I,s)$ is a Recoverable Support. 
The y-axis displays the proportion of Recoverable Supports versus all tested pairs.
For Mercedes Benz-frames only results for $n$ odd are displayed.
The formula \eqref{Eq:EstRandVertices} is plotted black, and formula \eqref{Eq:BoundMB} is displayed as the red line.
The lower bound from Corollary \ref{Co:LowerBound} is displayed in dashed red.}
\end{figure}

For the rest of this section we regard \textit{equiangular tight frames} $\{a_i\}_{1\le i\le n}$ in $\mathbb{R}^m$, where the vector $a_i$ forms the $i$-th column of the $m\times n$-matrix.
Among other things, these frames have the property that any pair of columns has the same inner product.
In case of minimally redundant matrices, i.e. $m=n-1$, the only equiangular tight frame is (up to rotation) the so-called \textit{Mercedes-Benz frame}, see \cite[Section 3.2]{MaPe09} and  \cite{IsPe06}.
Particular Mercedes-Benz frames have an additionally property: 
Each row of such a matrix has the mean value equal to zero, in other words, the kernel is spanned by the vector of all ones.
This property can be used to give the exact number of Maximal Recoverable Supports.
Let $n$ be odd.
Since any $v\in\rg{A^T}$ has the mean value zero, any vertex of $P$ has the same property.
We construct these vertices combinatorically by choosing an index set $J\subset\{1,...,n\}$ with $|J| = (n-1)/2$; there are $\binom{n}{(n-1)/2}$ different possibilities choosing $J$. 
Further there are $(n+1)/2$ different possibilities choosing one $l\in\{1,...,n\}\backslash J$. 
For, say, the Mercedes-Benz frame $A\in\mathbb{R}^{n-1\times n}$ it holds that $v\in\mathbb{R}^n$, with $v_i = 1$ for $i\in J$ and $v_l = 0$ as well as the remaining entries having the value $-1$, is an vertex of $P$. 
Hence 
\begin{align}
  \Lambda(A,n-1) = \left(\frac{n+1}2\right)\binom{n}{\frac{n-1}2}.\label{Eq:BoundMB}
\end{align}
Using the same argument for $n$ even, we get $\Lambda(A,n-1) = 0$ but $\Lambda(A,n-2) = (n/2)\binom{n}{n/2}$.
Keeping in mind that the combinatorical amount increases with a decreasing number of $\pm1$, we can construct any Recoverable Support of $A$ with any size, e.g. for $n$ even it holds 
\[\Lambda(A,n-2) = \left(\frac{n-1}2\right)\left(\frac{n+1}2\right)\binom{n}{\frac{n-1}2}.\]
The theoretical results so far are illustrated in Figure \ref{Fig:TheoResults} by Monte Carlo experiments with Mercedes-Benz frames and randomly drawn Gaussian matrices.
 % compared to the lower bound from Corollary \ref{Co:LowerBound}, the estimated value in \eqref{Eq:EstRandVertices} and the value \eqref{Eq:BoundMB} for Mercedes-Benz frame related to the empirical probability.
One may observe that the empirical results agree with the theorectical statements.

For $n$ even we can also construct a matrix $A\in\mathbb{R}^{n-1\times n}$ similar to the formula \eqref{Eq:BoundMB}, this will be revisited in Section \ref{Se3:Bounds}.
\begin{Lemma}\label{Le:ConstNEven}
 Let $A\in\mathbb{R}^{m\times n}$ then there exists a matrix $B\in\mathbb{R}^{m+1\times n+1}$ such that $\Lambda(B,m+1) = 2\Lambda(A,m)$.
\end{Lemma}

\begin{Proof}
Consider the set $W = \{w : w\mbox{ satisfies }\eqref{Conditions:l1}\mbox{ for some }(I,s)\}$ and for $\alpha\neq 0$ the matrix 
\[B = \left[\begin{array}{cc}A & 0\\0 & \alpha\end{array}\right].\]
Then for any $w\in W$ the elements $w^{(1)}=(w,\alpha^{-1})^T, w^{(2)}=(-w,\alpha^{-1})^T$ satisfy \eqref{Conditions:l1} for $B$.
Hence there are $2\Lambda(A,m)$ Recoverable Supports of $B$ with size $m+1$.
\end{Proof}

Since for $n$ even it holds
\[\left(n+2-\frac n2\right)\binom{n+2}{\frac n2} = 2\left(n+1-\frac n2\right)\binom{n+1}{\frac n2},\]  
and by denoting $\lfloor\cdot\rfloor$ as the \textit{Floor function}, we can state matrices $A\in\mathbb{R}^{n-1\times n}$ satisfying
\[\Lambda(A,n-1) = \left(n-\left\lfloor\frac n2\right\rfloor\right)\binom{n}{\lfloor\frac n2\rfloor}.\]
This formula will be important in Corollary \ref{Co:ONeil}.

Up to here, the partial order in the set of all Recoverable Supports of a certain matrix has not been used.
The following lemma enters this subject. 
It will be helpful for bounding $\Lambda$ and $\Xi$ and further gives some characteristics about the actual recoverability which is the number of Recoverable Supports in proportion to the total number of $(n-k)$-faces of $C^n$ (where $k$ is the size of the appropriate Recoverable Support).
\begin{Lemma}\label{Le:RatioIncomOutgoVert}
Let $A\in\mathbb{R}^{m\times n}$, then for any $k\le\mathrm{rank}(A)$ with $\Lambda(A,k)\neq 0$, there exists a positive number $\lambda\le2(n-k+1)$ satisfying
\[\lambda\Lambda(A,k-1) = k\Lambda(A,k).\]
\end{Lemma}
\begin{Proof}
Regarding the lattice of all Recoverable Supports of $A$, Theorem \ref{Th:poset} states that any Recoverable Support with size $k$ is adjacent to $k$ Recoverable Supports with size $k-1$, i.e. the number $k\Lambda(A,k)$ states the number of all adjacences between Recoverable Supports with size $k$ and $k-1$.
Hence, there is a positive number $\lambda$ satisfying the desired equation.
Any Recoverable Support $(I,s)$ with size $k-1$ is adjacent to no more than $2(n-k+1)$ Recoverable Supports with size $k$, since $|I^c| = n - k + 1$ and each new $s_j, j\in I^c$, in a Recoverable Support with size $k$ can adopt both signs: a positve or a negative sign.
Hence, it holds $\lambda\le2(n-k+1)$.
\end{Proof}

The number $\lambda$ from Lemma \ref{Le:RatioIncomOutgoVert} states the averaged number of outgoing adjacences from a Recoverable Support with size $k-1$ to Recoverable Supports with size $k$.
The upper bound for $\lambda$ implies a statement for the probability that an appropriate pair $(I,s)$ is a Recoverable Support.
\begin{Proposition}\label{Pr:Recov}
Let $A\in\mathbb{R}^{m\times n}$, then the mapping 
\begin{align}
 k\mapsto\left[2^k\binom{n}{k}\right]^{-1}\Lambda(A,k)\label{Map:Recov}
\end{align}
is monotonically nonincreasing.
\end{Proposition}

\begin{Proof}
 Assume there is $k\le\mbox{rank}(A)$ satisfying
\[\left[2^{k-1}\binom{n}{k-1}\right]^{-1}\Lambda(A,k-1)>\left[2^k\binom{n}{k}\right]^{-1}\Lambda(A,k).\]
Since there is $\lambda\in\mathbb{R}$ such that $\lambda\Lambda(A,k-1) = k\Lambda(A,k)$, it holds $\lambda > 2(n-k+1)$, which is a contradiction to Lemma \ref{Le:RatioIncomOutgoVert}.
\end{Proof}

The mapping \eqref{Map:Recov} states the ratio between the actual number of Recoverable Supports of $A$ with size $k$ and the total number of all pairs $(I,s)$ with $I\subset\{1,...,n\}, s\in\{-1,+1\}^I$, previously introduced as the recoverability.
The second proposition aims at an actual number of $\Lambda$ for sparsity $\mbox{rank}(A)-1$ if the number of Maximal Recoverable Supports is known. 
\begin{Proposition}\label{Pr:SparsityK-1}
 Let $A\in\mathbb{R}^{m\times n}$ with rank $l$ and assume $\Lambda(A,l)\neq0$.
Then $\Lambda(A,l-1) = \frac l2\Lambda(A,l)$.
\end{Proposition}

\begin{Proof}
Regarding any Recoverable Support $(I,s)$ with size $l-1$, it holds that the null space of $A_I^T$ is one-dimensional.
Since there exist at least one Recoverable Support with size $l$, we can enlarge it, due to Theorem \ref{Th:poset}, in two different directions.
\end{Proof}

Proposition \ref{Pr:SparsityK-1} states another interesting fact about the number of Maximal Recoverable Support:
Noticing that all values of $\Lambda$ are even due to the symmetry of the underlying polytope, we observe for an odd rank $l$ of a matrix $A$ that $\Lambda(A,l)$ is divisible by four or even a higher even number.

\subsection{Bounds and Values for $\Xi$}\label{Se3:Bounds}
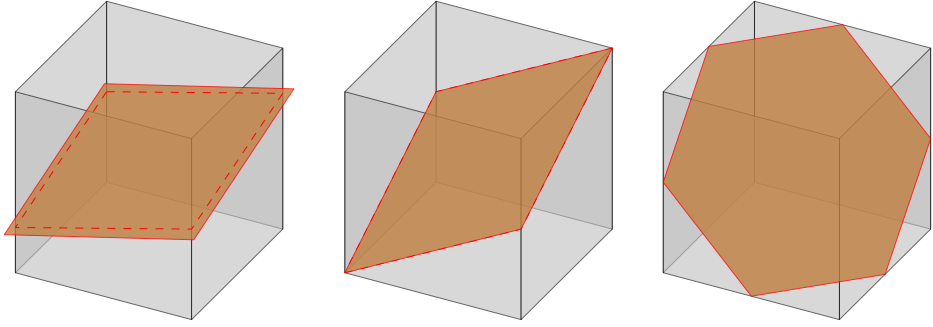
\begin{figure}
 \begin{tikzpicture}[ x = {(-0.5cm,-0.5cm)},
			y = {(0.9659cm,-0.25882cm)},
			z = {(0cm, 1cm)},
			scale = 1.2,
			color = {lightgray}]

\tikzset{facestyle/.style={fill=lightgray, draw=black,opacity=.6, very thin, line join = round}}

% Begin Cube
\begin{scope}[canvas is zy plane at x=0]
	\path[facestyle] (0,0) rectangle (2,2);
\end{scope}

\begin{scope}[canvas is zx plane at y=0]
	\path[facestyle] (0,0) rectangle (2,2);
\end{scope}

\begin{scope}[canvas is zy plane at x=2]
	\path[facestyle] (0,0) rectangle (2,2);
\end{scope}

\begin{scope}[canvas is zx plane at y=2]
	\path[facestyle] (0,0) rectangle (2,2);
\end{scope}
% End Cube

% Begin Surface
\draw[fill=brown,draw=red,opacity=.8,very thin,line join=round]
(2.05,-0.1,0.42) --
(2.05,2.05,0.92) --
(-0.05,2.1,1.55) --
(-0.05,-0.05,1.05) --cycle
;
% End Surface

% 
% Begin Intersection Surface and Cube
\draw[very thin,red,dashed,line join=round]
(0,0,1) -- 
(2,0,0.5) --
 (2,2,1) --
 (0,2,1.5) --cycle
;
% End Intersection
\end{tikzpicture}
\hfill
\begin{tikzpicture}[ x = {(-0.5cm,-0.5cm)},
			y = {(0.9659cm,-0.25882cm)},
			z = {(0cm, 1cm)},
			scale = 1.2,
			color = {lightgray}]

\tikzset{facestyle/.style={fill=lightgray, draw=black,opacity=.6, very thin, line join = round}}

% Begin Cube
\begin{scope}[canvas is zy plane at x=0]
	\path[facestyle] (0,0) rectangle (2,2);
\end{scope}

\begin{scope}[canvas is zx plane at y=0]
	\path[facestyle] (0,0) rectangle (2,2);
\end{scope}

\begin{scope}[canvas is zy plane at x=2]
	\path[facestyle] (0,0) rectangle (2,2);
\end{scope}

\begin{scope}[canvas is zx plane at y=2]
	\path[facestyle] (0,0) rectangle (2,2);
\end{scope}
% End Cube

\draw[very thin,red,dashed,line join=round]
(2,0,0) --
 (2,2,1) --
 (0,2,2) --
 (0,0,1)--cycle;

\draw[fill=brown,draw=red,opacity=.8,very thin,line join=round]
(2,0,0) --
 (2,2,1) --
 (0,2,2) --
 (0,0,1)--cycle
;
\end{tikzpicture}
\hfill
\begin{tikzpicture}[ x = {(-0.5cm,-0.5cm)},
			y = {(0.9659cm,-0.25882cm)},
			z = {(0cm, 1cm)},
			scale = 1.2,
			color = {lightgray}]

\tikzset{facestyle/.style={fill=lightgray, draw=black,opacity=.6, very thin, line join = round}}

% Begin Cube
\begin{scope}[canvas is zy plane at x=0]
	\path[facestyle] (0,0) rectangle (2,2);
\end{scope}

\begin{scope}[canvas is zx plane at y=0]
	\path[facestyle] (0,0) rectangle (2,2);
\end{scope}

\begin{scope}[canvas is zy plane at x=2]
	\path[facestyle] (0,0) rectangle (2,2);
\end{scope}

\begin{scope}[canvas is zx plane at y=2]
	\path[facestyle] (0,0) rectangle (2,2);
\end{scope}
% End Cube

\draw[fill=brown,draw=red,opacity=.8,very thin,line join=round]
(1,0,2) --
(2,0,1) -- 
(2,1,0) --
(1,2,0) --
(0,2,1) --
(0,1,2) -- cycle
;
\end{tikzpicture}
\caption{Examples for Sections $\rg{A^T}\cap C^3$. \textit{Left:} Regular Section; \textit{Center:} Not a Regular Section; \textit{Right:} Regular Secction with Mercedes Benz Frame}
\label{Fi:ExSection}
\end{figure}

In this subsection, we give bounds and values for the largest possible number of Recoverable Supports of all matrices with a certain size, i.e.~$\Xi$ (cf. Definition~\ref{Def:Bounds}).

It is obvious that we can slice the three dimensional cube $C^3$ with a hyperplane in maximal six edges, see Figure~\ref{Fi:ExSection}.
As Figure \ref{Fi:ExSection} prompts it is not possible to slice less than four edges without failing the origin, the graphics in the middle shows that it is possible to touch also vertices of the hypercube. 
Despite Theorem \ref{Th:GeometricalInterpretL1} implies that the results from the field \textit{cross-sections of a hypercube} can be used for our issues, these results often require a regular cross-section while, in general, the section $\rg{A^T}\cap C^n$ is not regular.
In contrast to lower bounds (cf. Section \ref{Se:ValLam}), results for an upper bound can be used, as regarded in the following of this subsection.
Note that McMullens Upper Bound Theorem \cite{M70} can not be used as a typical choice, since it exceeds the trivial bound.

Firstly we give an upper bound for $\Xi$ if $k$ is large.
This result is already known \cite[Corollary 1.3]{Do04} (cf. Corollary \ref{Co:NecK}) in the field of \textit{Sparse Reconstruction}.
\begin{Corollary}
Let $0<m<n-1$. If $k>\frac{m+1}3$ then $\Xi(m,n,k)<2^k\binom{n}{k}$.
\end{Corollary}
\begin{Proof}
 This result follows from \cite{La78,MuSh68} and the second statement in Theorem \ref{Th:GeometricalInterpretL1}.
\end{Proof}

Considering minimally redundant matrices, remind $m=n-1$, we get the following value for Maximal Recoverable Supports.
\begin{Corollary}\label{Co:ONeil}
It holds
\[\Xi(n-1,n,n-1) = \left(n-\left\lfloor\frac n2\right\rfloor\right)\binom{n}{\lfloor\frac n2\rfloor}.\]
\end{Corollary}%
\begin{Proof}
 Combining \cite{Ne71} and Statement 2 of Theorem \ref{Th:GeometricalInterpretL1} proves the result.
\end{Proof}

In Section \ref{Se:ValLam} we have seen that the Mercedes-Benz frame with an odd number of columns and the construction in Lemma \ref{Le:ConstNEven} reaches this value.
Additionally, with the \textit{mutual coherence} slightly more than half of the values $\Xi(n-1,n,k)$ for variable $k$ are known from the following result.
\begin{Corollary}\label{Co:MuCo}
It holds 
\[\Xi(m,n,k) = 2^k\binom{n}{k}\mbox{ if }k<\frac 12\left(1+\sqrt{\frac{m(n-1)}{n-m}}\right).\]
\end{Corollary}%
\begin{Proof}
  This follows from \cite{DoHu01,StHe03}.
\end{Proof}

The bound in Corollary \ref{Co:MuCo} can be reached by equiangular tight frames, see \cite{StHe03}.
As a further consequence of the bound in Lemma \ref{Le:RatioIncomOutgoVert}, the following proposition delivers an upper bound for $\Xi$.

\begin{Proposition}
 For $k\le m$ it holds 
\[\Xi(m,n,k)\le\frac {2(n-k+1)}k\Xi(m,n,k-1).\]
\end{Proposition}

\begin{Proof}
 Assume there is $k\le m$ for $A\in\mathbb{R}^{m\times n}$ with $\Xi(m,n,k-1)=\Lambda(A,k-1)$ and $\tilde A\in\mathbb{R}^{m\times n}$ with $\Xi(m,n,k)=\Lambda(\tilde A,k)$ satisfying
\[\Lambda(\tilde A,k) > \frac {2(n-k+1)}k\Lambda(A,k-1),\]
then it holds
\[\frac {2(n-k+1)}k\Lambda(A,k-1)<\Lambda(\tilde A,k)\le\frac {2(n-k+1)}k\Lambda(\tilde A,k-1)\]
with Lemma \ref{Le:RatioIncomOutgoVert}, which is a contradiction to $\Xi(m,n,k-1)=\Lambda(A,k)$.
\end{Proof}

Similarly to the value $\Lambda$, the latter result implies further statements about $\Xi$, which are similar to Propositions \ref{Pr:Recov} and \ref{Pr:SparsityK-1}.
\begin{Corollary}\label{Co:Xil-1} It holds $\Xi(m,n,m-1) = \frac m2\Xi(m,n,m)$.
\end{Corollary}

Additionally, we get a similar statement to Proposition \ref{Pr:Recov} about an upper bound of the recoverability.
\begin{Corollary}\label{Co:RecCurve_Xi}
The mapping
\[k\mapsto\left[2^k\binom{n}{k}\right]^{-1}\Xi(m,n,k)\]
is monotonically nonincreasing.
\end{Corollary}

To the end of this section, we develop a heuristic upper bound of $\Xi$.
Considering $\lambda$ in Lemma~\ref{Le:RatioIncomOutgoVert}, we can establish an upper bound of $\Xi$ by assuming that $\lambda$ can be bounded from below, i.e. $\lambda\ge2(l-k+1)$ for matrices with rank $l$.
Conveniently, we derive this heuristic bound for full rank, minimally redundant matrices $A$, i.e. $l=n-1$, but the construction can be adapted straightforward to other instances.
Assume $\lambda\ge2(n-k)$, then for a positive integer $v < n$ it follows
\[\Lambda(A,n-1) \ge 2^{v-1}\frac{(v-1)!(n-v)!}{(n-1)!}\Lambda(A,n-v)\]
by applying the lower bound recursively.
Through substituting $k= n-v$ and bounding $\Lambda(A,k-1)$ by Corollary~\ref{Co:ONeil}, we obtain
\[\Lambda(A,k) \le 2^{k+1-n}\binom{n-1}k\left(n-\left\lfloor\frac n2\right\rfloor\right)\binom n{\left\lfloor\frac n2\right\rfloor}.\]
Since the right-hand side of the latter inequality exceeds the trivial bound $2^k\binom nk$ for small $k$, we postulate the following heuristic upper bound:
\begin{align}
\Xi(n-1,n,k)\le\min\left\{2^k\binom{n}{k}, 2^{k+1-n}\binom{n-1}{k}\left(n-\left\lfloor\frac n2\right\rfloor\right)\binom{n}{\left\lfloor\frac n2\right\rfloor}\right\}.\label{InEq:UpBoundBeta}
\end{align}
In general, the inequality $\lambda\ge2(l-k+1)$ is not true, but we motivate this bound by the observation that the transition from all pairs $(I,s)$ are Recoverable Supports to none of the pairs $(I,s)$ are Recoverable Supports is rapid, e.g. \cite{TsDo06,BrDoEl09}, and, furthermore, this bound is true and strict in case that $k=l$, cf. Corollary \ref{Co:Xil-1}.
As far as we know, there is no matrix exceeding this heuristic; it will be considered in the computational experiments in Section 5.
Moreover, this bound is also strict due to Corollary~\ref{Co:ONeil}, \ref{Co:MuCo} and \ref{Co:Xil-1} for some values of $k$.

In the context of the Hasse Diagram of all Recoverable Supports, the maximum $\Xi$ also states a geometrical question: what is the maximal $\lambda$ such that the ratio $\lambda/k$ between the outgoing edges of all Recoverable Supports with size $k-1$ and the incoming edges of all Recoverable Supports with size $k$?
Results for this question would give further insights about $\Xi$ and improve a non-trivial upper bound.
Combining Corollary~\ref{Co:RecCurve_Xi} with Corollary \ref{Co:MuCo} delivers an interesting insight for $A\in\mathbb{R}^{m\times n}$: if it holds $\Lambda(A,\tilde k) = 2^k\binom{n}{\tilde k}$ for some $\tilde k\le m$, then equality holds in Lemma~\ref{Le:RatioIncomOutgoVert} for all $k\le\tilde k$, and also for $\Xi$.

%%% Local Variables: 
%%% mode: latex
%%% TeX-master: "main"
%%% End: 

%% file: section4_algorithm.tex
\section{Computing a Recoverable Support}
\label{sec:comp-rec-supp}
In general, generating test instances for computational experiments is an expensive problem in Basis Pursuit. 
Even for, say, Gaussian matrices, where one only has to find an instance satisfying the optimality condition for $\ell_1$ minimization derived by its subdifferential, it is not straightforward to find a suitable $x^*$ satisfying \eqref{l1} if the desired $x^*$ shall not be very sparse.

One na\"ive way to generate a test instance is to choose an arbitrary $k$-sparse vector $x^*$, solve~(\ref{l1}) with some solver and then check whether the solution is equal to $x^*$.
This may work well for small $k$ but usually becomes computationally expensive for larger $k$. Moreover, this construction suffers from a ``trusted method bias'', i.e. the method used to solve~(\ref{l1}) may work better on instances which inherit a particular structure (something which may not be under control of the experimenter).
Another approach has been proposed in~\cite{Lo12}: Choose a pair $(I,s)$ randomly and construct a dual certificate, i.e. find $w$ as in~(\ref{Conditions:l1}).
This problem could be seen as a convex feasibility problem~\cite{bauschke1996convexfeasibility} and can be solved, e.g., by alternating projections as outlined in~\cite{Lo12}.
This approach often leads to dual certificates $w$ such that the value $\|A_{I^c}^Tw\|_\infty$ is close to one and hence, the result may not be trustworthy due to numerical errors. A more favorable way to check the reconstructability using~(\ref{Conditions:l1}) would be to check if for some $(I,s)$ the optimal value of 
\begin{align}
 \min_w\|A_{I^c}^Tw\|_\infty \mbox{ subject to }A_I^Tw=s_I\label{MinProg:DualCert}
\end{align}
is less or equal one. Similar to the $\ell_1$ minimization problem~(\ref{l1}), this may be cast as a linear program. However, there are import differences to the na\"ive approach: First, the number of variables is $m$ which may be much smaller than $n$. Moreover, one does not rely on the entries of $x^*$ but only on its sign and the support.

However, in all the above methods one generates some trial support $(I,s)$ and then checks whether it is recoverable.
Derived from Corollary \ref{Co:EstRandVertices}, the probability for an appropriate pair $(I,s), |I|=n-1,$ being a Maximal Recoverable Support of a randomly drawn Gaussian matrix of the size $(n-1)\times n$ tends to zero for huge $n$.
Hence, one may never find any $(n-1)$-sparse vector by any trial-and-error method and a similar conclusion is true for $k$-sparse vectors for $m\times n$ matrices if $k$ is sufficiently large.
But in view of Theorem~\ref{Th:poset}, there is a systematic way to generate Recoverable Supports $(I,s)$ with maximal size by selecting a $1$-sparse recoverable vector, computing a corresponding dual certificate and incrementally increasing the support while maintaining a valid dual certificate (according to Theorem~\ref{Th:poset}, 1.).
The method is outlined in Algorithm \ref{Algo:ComputingRS}.
Note that there is considerable freedom in lines \ref{Algo1:ChooseBasisVector} and \ref{Algo1:MinProb} of the algorithm on how to continue.

\begin{algorithm}[htb]
\DontPrintSemicolon
\SetAlgoLined
\SetKwInOut{Input}{Input}\SetKwInOut{Output}{Output}
\SetKwData{Jj}{$|J|\le k$ and $A_J$ has full rank}
 \Input{$A\in\mathbb{R}^{m\times n}, k\le\mbox{rank}(A)$}
 \Output{Recoverable Support $(I,s)$ of $A$ with size $k$}
 \BlankLine
 \emph{$a_k = \arg\max_{a_i}\|a_i\|_2^2$}\tcp*[f]{The $i$-th column of $A$ is denoted by $a_i$}\;
 \emph{$w\leftarrow\|a_k\|_2^{-2}a_k$}\;
 \emph{$s\leftarrow A^Tw$}\label{Algo1:ReSuSize1}\;
 \emph{$I\leftarrow\{k\}$}\label{Algo1:ReSuSize1a}\;
 \emph{$I^c\leftarrow\{1,...,n\}\backslash\{k\}$}\label{Algo1:ReSuSize1b}\;
 \While{$|I|<k$}{\label{Algo1:WhileLoop}
   \emph{Choose a vector $y\in\ker{A_I^T}$}\label{Algo1:ChooseBasisVector}\;

   \emph{Choose $\lambda\in\mathbb{R}$ such that $\|A_{I^c}^T(w+\lambda y)\|_\infty = 1$}\label{Algo1:MinProb}\;
   \emph{$J\leftarrow\{i : |a_i^T(w+\lambda z)| = 1\}$}\;
   \uIf{\Jj}{\label{Algo1:IfCause}
     \emph{$I\leftarrow J$}\;
     \emph{$I^c\leftarrow\{1,...,n\}\backslash I$}\;
     \emph{$w\leftarrow w+\lambda y$}\;
     \emph{$s\leftarrow A^Tw$}\label{Algo1:WalkTheSubset}\;
   }
   \lElse{
     \emph{Return to line \ref{Algo1:ChooseBasisVector}}\label{Algo1:ReturnToChooseBasisVector}\;
   }
 }
 \caption{Computing a Recoverable Support}\label{Algo:ComputingRS}
\end{algorithm}

Algorithm \ref{Algo:ComputingRS} is designed for arbitrary matrices of arbitrary sizes.
However, it is possible that the algorithm does not deliver a desired Recoverable Support if it gets stuck in line~\ref{Algo1:IfCause}.
To protect against these cases the method could be extended by including the second statement of Theorem \ref{Th:poset};
this extension would deliver more freedom to jump between different index sets $I$ but requires elaborate bookkeeping of previously visited index sets.
We experienced that this extension is not necessary in most cases.

The first issue about the algorithm might be the question, for what kind of matrices does the method compute a Recoverable Support. 
Theorem \ref{Th:Existence} gives an answer: matrices which columns with maximal Euclidean norm are pairwise linearly independant.
The construction in the proof of Theorem \ref{Th:Existence} for a Recoverable Support with size one is used in the first three lines.
Hence, for these  matrices the variable $s$ in line \ref{Algo1:ReSuSize1} has only one entry equal to one in absolute value, the rest of the absolute entries are less than one; this occasions the clauses in line \ref{Algo1:ReSuSize1a} and \ref{Algo1:ReSuSize1b}.
%In general, an extension to other startings is not necessary since Theorem \ref{Th:Existence} covers a huge number of types of matrices.

Theorem \ref{Th:GeometricalInterpretL1} gives a geometrical interpretation of Algorithm \ref{Algo:ComputingRS}. 
In line \ref{Algo1:ReSuSize1} we start on one facet of the hypercube and by line \ref{Algo1:WalkTheSubset} we walk along the range of the transposed matrix to the next lower-dimensional face of the hypercube.
Consequently, the method requires at least $k-1$ iterations for computing a Recoverable Support with size $k$.
Experiences show that mostly only $k-1$ iterations are required.
The if-clause in line \ref{Algo1:IfCause} saves for being \textit{stuck} in an unsuitable face.

In any iteration step of the while loop, an element of the corresponding null space is chosen.
To choose such a vector it is advantageous to maintain an orthonormal basis for the kernel of $A_I^T$ during the iteration in the form of some decomposition.
In our setting, we are calling up a rank one update to a QR decomposition.
In the worst-case scenario it may happen that one needs to check several vectors $y$ in line~\ref{Algo1:ChooseBasisVector}, however, using an orthonormal basis of the kernel one can just try all of the basis vectors one after another.
This worst case would lead to an iteration number $\mathcal{O}(l^2)$ for computing a Recoverable Support with size $l$.
Actually, we were not able to construct such an instance and usually the iteration number is $\mathcal{O}(l)$.
Our setting of this method, implemented as a MATLAB program, can be found online at \homepage.

% For many Recoverbable Support $(I,s)$ there exists many elements $w$ satisfying \eqref{Conditions:l1}, hence one can also extend Algorithm \ref{Algo:ComputingRS} by adding a certain condition, e.g. stated in \eqref{MinProg:DualCert} for $s=\sign{x^*_I}$, for regarding a specific values of $w$ after the while-loop starting in line \ref{Algo1:WhileLoop}.

% This may be interesting for exploring characteristics for different dual certificates.

%%% Local Variables: 
%%% mode: latex
%%% TeX-master: "main"
%%% End: 

%% file: section5_computational_experiments.tex
\section{Computational Experiments}
In this section, we present computational experiments for the topics of the previous sections.
The optimization problem \eqref{MinProg:DualCert} delivers an alternative method to perform numercial experiments in Basis Pursuit.
A comparison of solving \eqref{MinProg:DualCert} and solving the $\ell_1$ minimization in \eqref{l1} will be done in the following subsection.
In Subsection \ref{Se:NRSCM} we will highlight the theorectical results from Section 3 with Monte Carlo experiments and will show the behaviour of the heuristic upper bound from \eqref{InEq:UpBoundBeta}.
% Further we will give empirical results on Algorithm \ref{Algo:ComputingRS} in Subsection \ref{Se:PerformAlgo}.
% In the last subsection, we extend \eqref{l1} to the analysis $\ell_1$ minimization problem and present results similar to Subsection \ref{Se:Com_l1_l_inf} and \ref{Se:NRSCM}.

All experiments were done with Matlab R2012b employed on a desktop computer with 4 CPUs, each Intel\textregistered\ Core\texttrademark\ i5-750 with 2.67GHz, and 5.8 GB RAM; the $\ell_1$ and $\ell_\infty$ minimization problems were solved as linear programs with Mosek 6.

In the Monte Carlo experiments it will be tested whether a pair $(I,s)$, with $I\subset\{1,...,n\}, s\in\{-1,1\}^I$, is a Recoverable Support of a given matrix.
The experiments were done as follows:
For a given matrix $A\in\mathbb{R}^{m\times n}$ and $k\le m$, we generate $I\subset\{1,...,n\}$ with $|I|=k$ randomly by choosing $I$ uniformly at random over $\{1,...,n\}$ and assure whether the submatrix $A_I$ has full rank through the Matlab function \texttt{rank}.
If $A_I$ has no full rank, then $(I,s)$ is not a Recoverable Support of $A$; otherwise we also choose $s\in\{-1,1\}^I$ randomly and solve the $\ell_\infty$ minimization problem \eqref{MinProg:DualCert} with $s = \sign{x^*_I}$.
If the optimization problem is feasible, solved with status 'optimal' and its optimization value is strictly less than one, the pair $(I,s)$ will be recorded as a Recoverable Support of $A$.
For each size $k$, we perform $M$ repetitions and average the results; the number $M$ varies from experiment to experiment and may be obtained from the descriptions to each experiment.
For reproducibility the code for all tests is at \homepage.

\subsection{Comparing $\ell_1$ and $\ell_\infty$ Solver in Mosek}\label{Se:Com_l1_l_inf}

To check whether a pair $(I,s)$ is a Recoverable Support, there are different methods, e.g. outlined in Section~\ref{sec:comp-rec-supp}.
In this subsection, we compare the na\"ive approach, i.e. solving \eqref{Conditions:l1} for some $x^*$ with the desired signum $s$, with solving \eqref{MinProg:DualCert}.
For comparision, we decided to perform a similar setup as in typical studies of the \textit{Phase Transition}, see e.g. \cite{DoTa09}.
We chose, as in \cite{DoTa09}, Gaussian matrices $A\in\mathbb{R}^{m\times n}$ for fixed $n=1600$ and varying $m$ such that $\delta=m/n\in(0,1]$ is chosen in forty equidistant steps.
The tests were realized as Monte Carlo experiments with varying $|I|=k$ such that for any $m$ the value $\rho=k/m\in(0,1]$ is chosen in forty equidistant steps.
For any triple $(m,n,k)$, we did the following testing.
We chose $A\in\mathbb{R}^{m\times n}$ as a randomly drawn Gaussian matrix, and performed the Monte Carlo sampling as described above by firstly check whether $(I,s)$ is a Recoverable Support of $A$, then choose $x^*$ with $\supp{(x^*)}=I, \sign{x^*}_I = s$, and solve Basis Pursuit with the right-hand side $Ax^*$.
This procedure is done with $M=10$ repetitions.
Remarkably, both approaches can be cast as solutions of linear programs and hence, we used the same solver for linear programs.
More precisely, testing whether $(I,s)$ is a Recoverable Support by solving \eqref{MinProg:DualCert} was implemented as a linear program and solved with the Mosek routine \texttt{mosekopt} with all tolerances set to default.
We decide that the pair $(I,s)$ is a Recoverable Support if $A_I$ has full rank, the optimization problem is feasible, it is solved with a status 'Optimal', and its objective value is is strictly less then $1-10^{-12}$.
On the other hand, we checked if $x^*$ satisfies \eqref{l1} by solving the constrained $\ell_1$ minimization as a linear program with the Mosek routine \texttt{mosekopt}; again all tolerances were set to default.
We judge a calculated solution $\tilde x$ to be exact if $\|\tilde x-x^*\|<10^{-5}$.

\begin{figure}
\begin{center}\includegraphics[width=60mm]{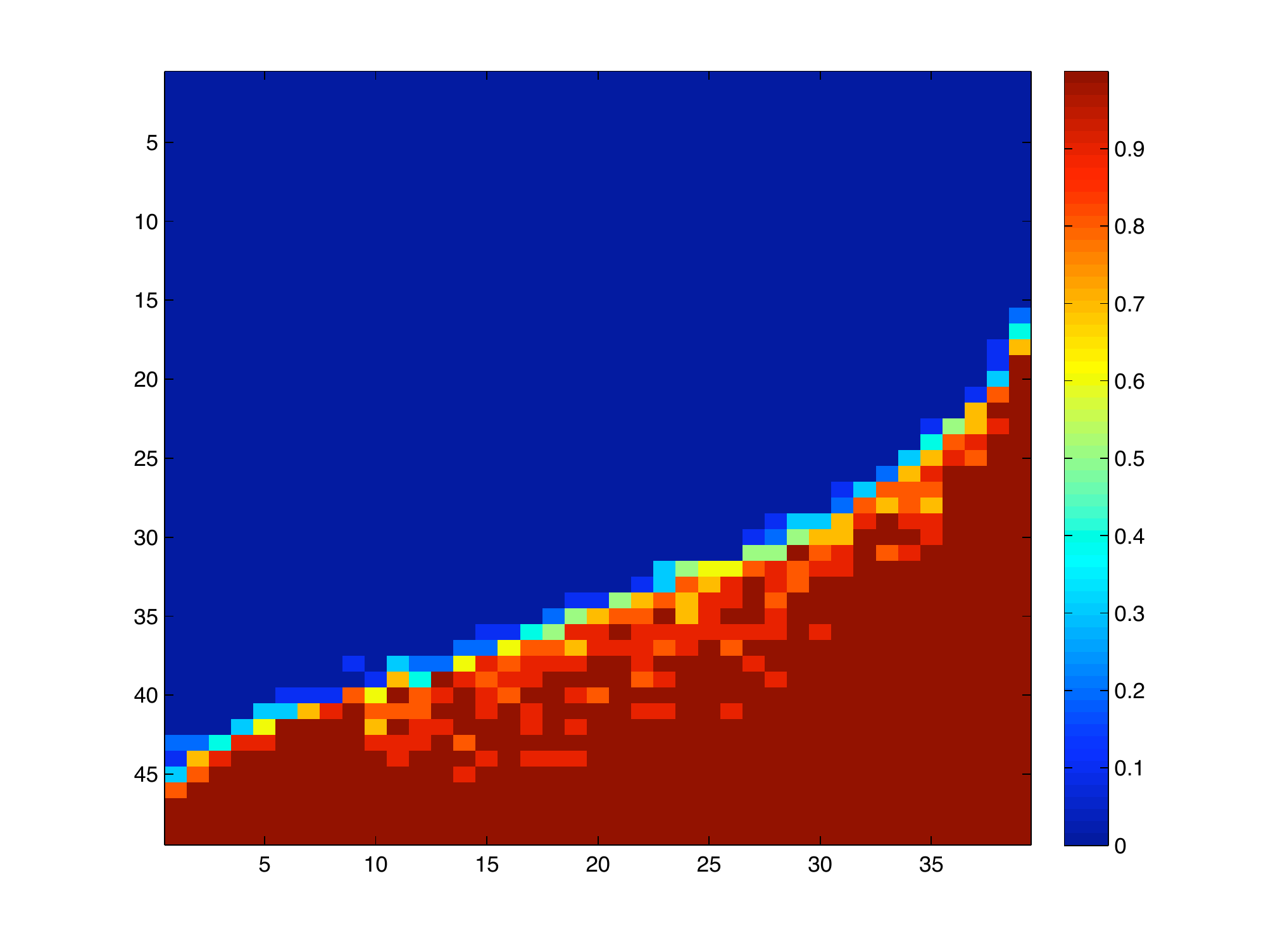}
\includegraphics[width=60mm]{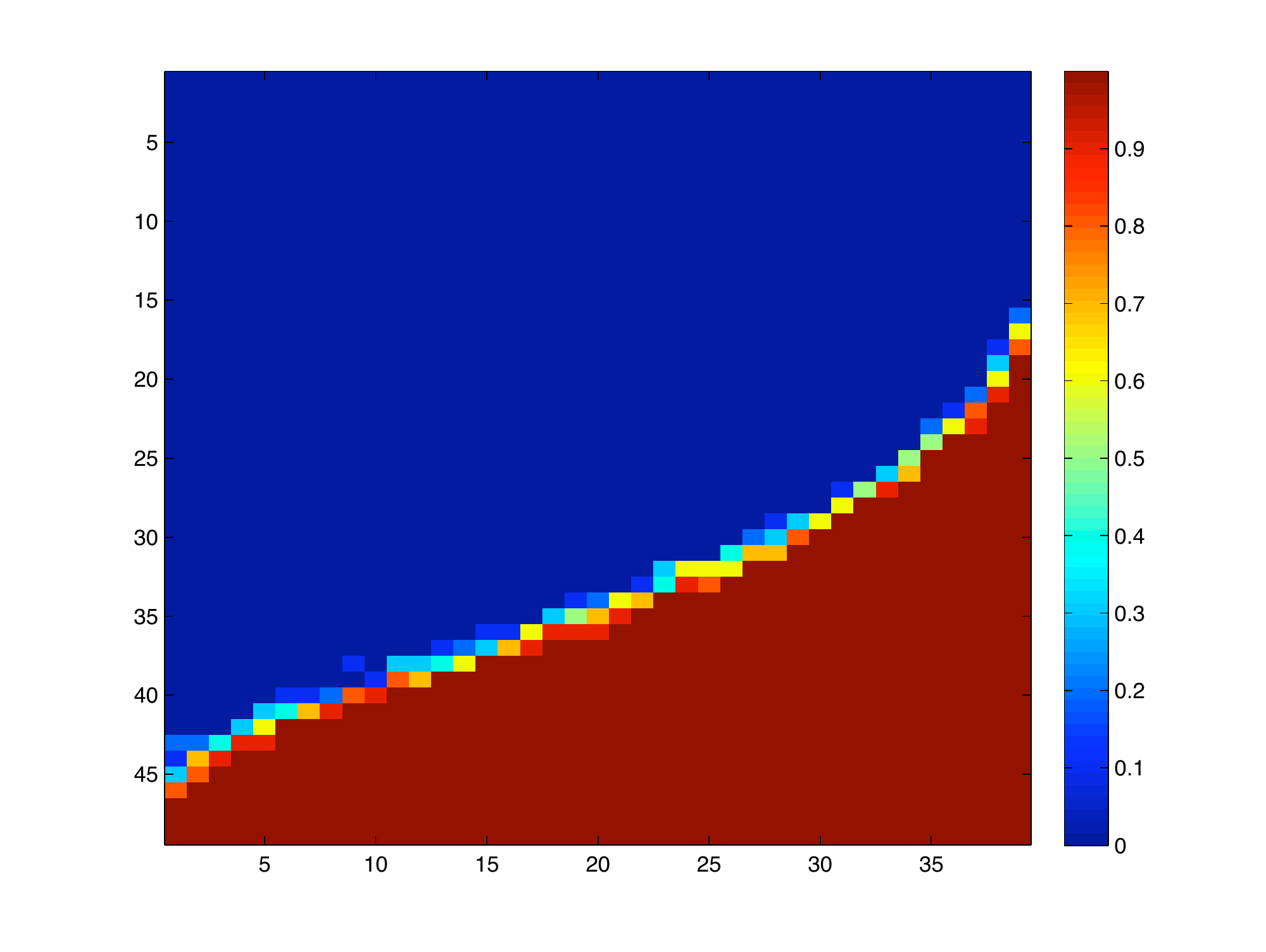}\end{center}
\caption{Averaged results from Monte Carlo experiments whether test instance is solved by $\ell_1$ minimization (left) and \eqref{MinProg:DualCert} (right). The values reach from zero (none of the instances where solutions) to one (all instances were solutions).}
\label{Fig:Time_Comp_PT}
\end{figure}

First we observe that all calculated solutions were solved with the status ``Optimal''.
Figure \ref{Fig:Time_Comp_PT} displays the averaged results of the decision whether a calculated solution of $\ell_1$ minimization is the desired solution (left) and a tested pair is a Recoverable Support (right).
The miss-fit between the figures comes from the fact that the solutions of~(\ref{l1}) are not accurate enough to fulfill the desired tolerance of $10^{-5}$. %, even though the tolerances for the solver were set quite high.
Relaxing the bound from $10^{-5}$ to $10^{-3}$ would lead to almost identical figures in this case but may lead to more errors in other circumstances.
Alternatively, instead of measuring the Euclidean distance between the calculated solution $\tilde x$ and the actual solution $x^*$, one may compare whether the support of $x^*$ and the support of $\tilde x$ coincide; however, to determine the support, another tolerance would be needed to identify the non-zero entries.
In perspective to previous experiments, e.g. \cite{DoTa09}, the results as in Figure~\ref{Fig:Time_Comp_PT} are as expected.
Further, we see agreement to previous testings as the phase transition between one to zero is displayed by the curve $\rho_F$ from Theorem \ref{Th:PropNeighbWeak} (cf. Figure \ref{Fig:PhaseTransition}).

For measuring the performance of both procedures, we measure the time it took to solve each linear program.
We excluded all operations to formulate the constraints of the linear programs from the time measurement.
Additionally before solving \eqref{MinProg:DualCert}, we checked whether $A_I$ has full rank and measure its duration.
If $A_I$ is not a full rank matrix, the problem \eqref{MinProg:DualCert} would not be solved.
Since we are only considering Gaussian random matrices, which are full spark matrices with probability 1, we could have skipped the testing of the rank (and we would have saved about 0.7 percent of the entire run time of the test) but we decided to present the test without any restrictions to specific test problems.

\begin{figure}
\begin{center}\includegraphics[width=100mm]{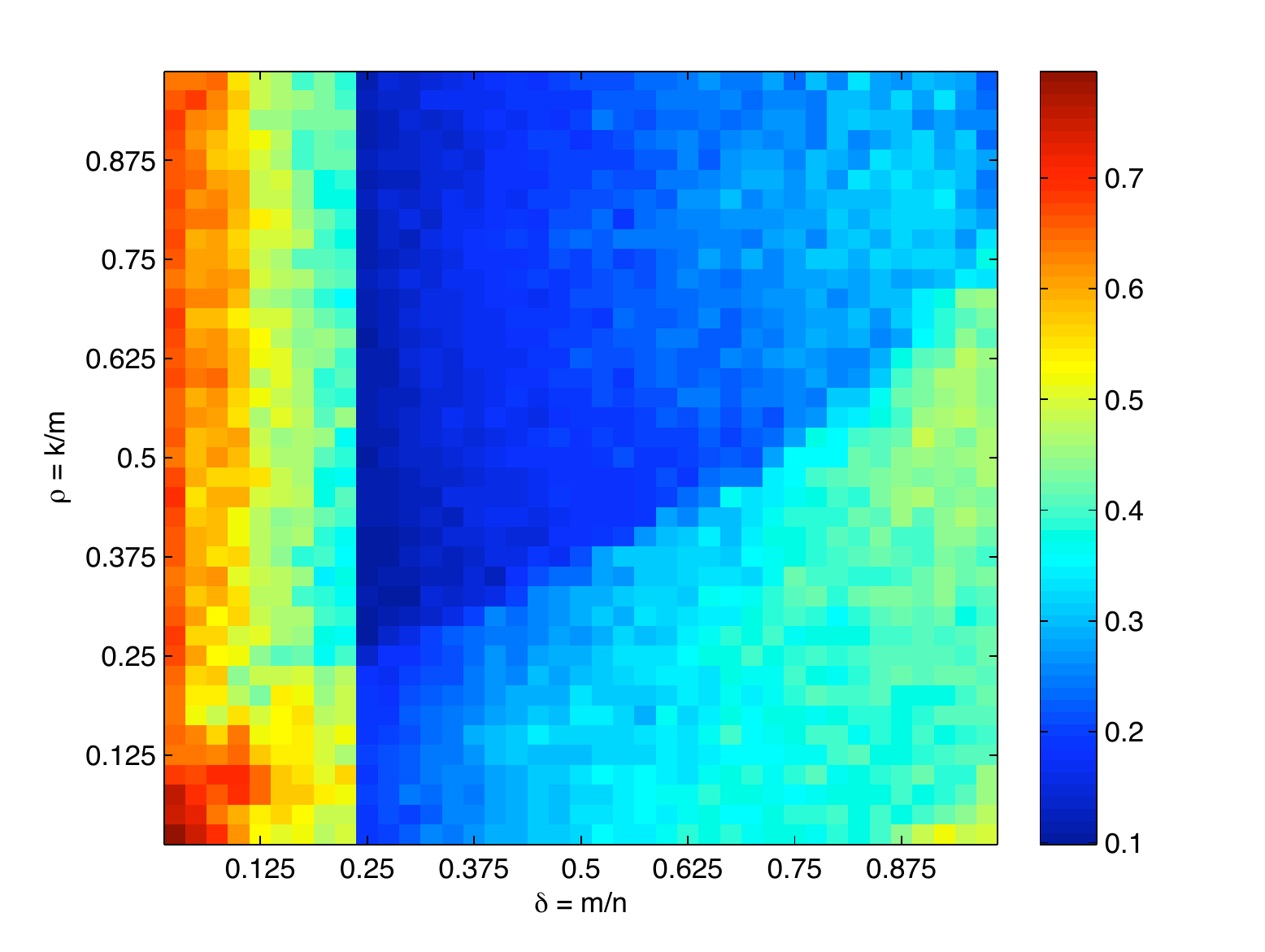}\end{center}
\caption{Comparision of duration performances between $\ell_1$ minimization and \eqref{MinProg:DualCert} and checking injectivity.}
\label{Fig:Time_Comp}
\end{figure}

\begin{figure}
\begin{center}\includegraphics[width=60mm]{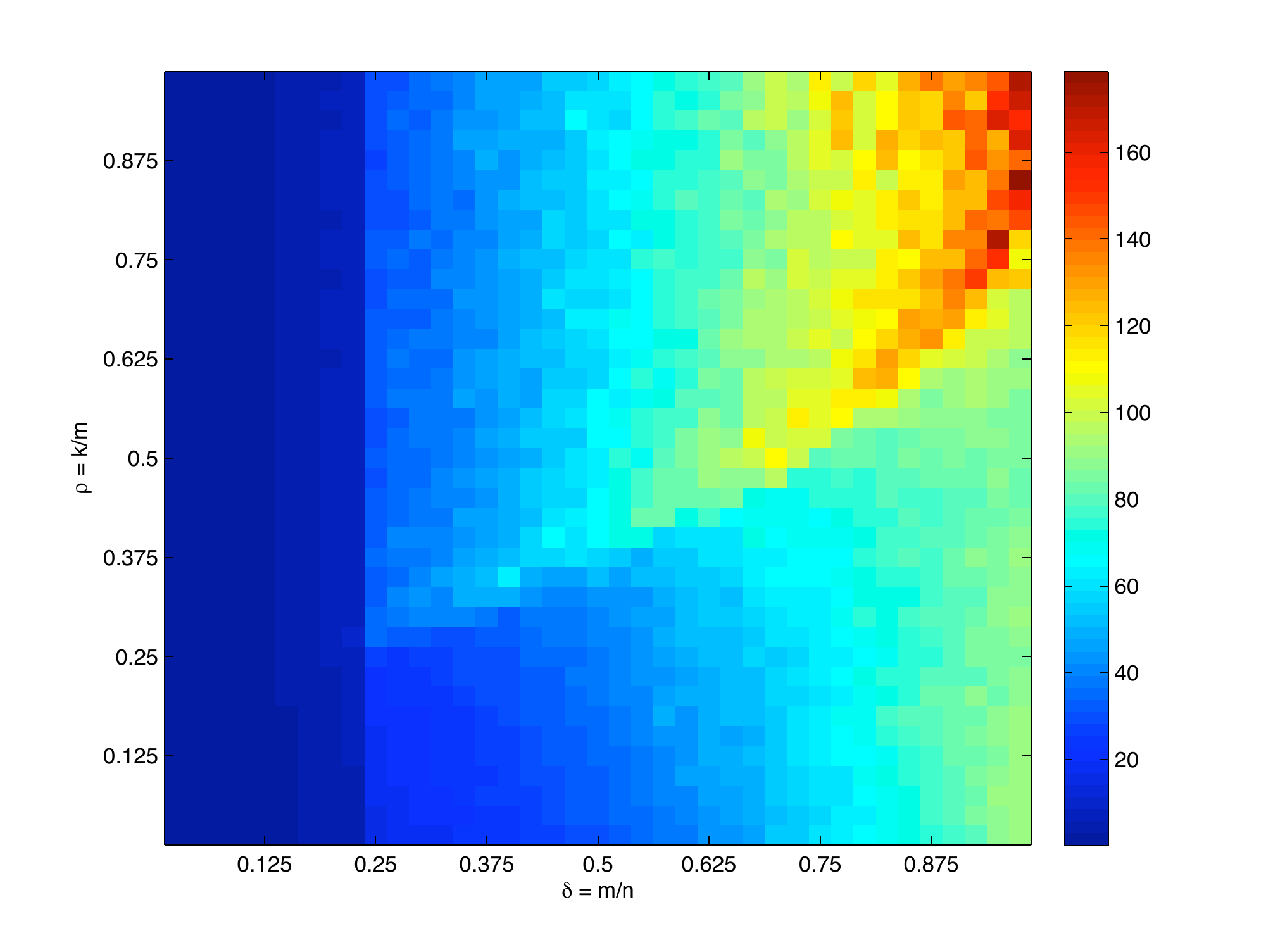}
\includegraphics[width=60mm]{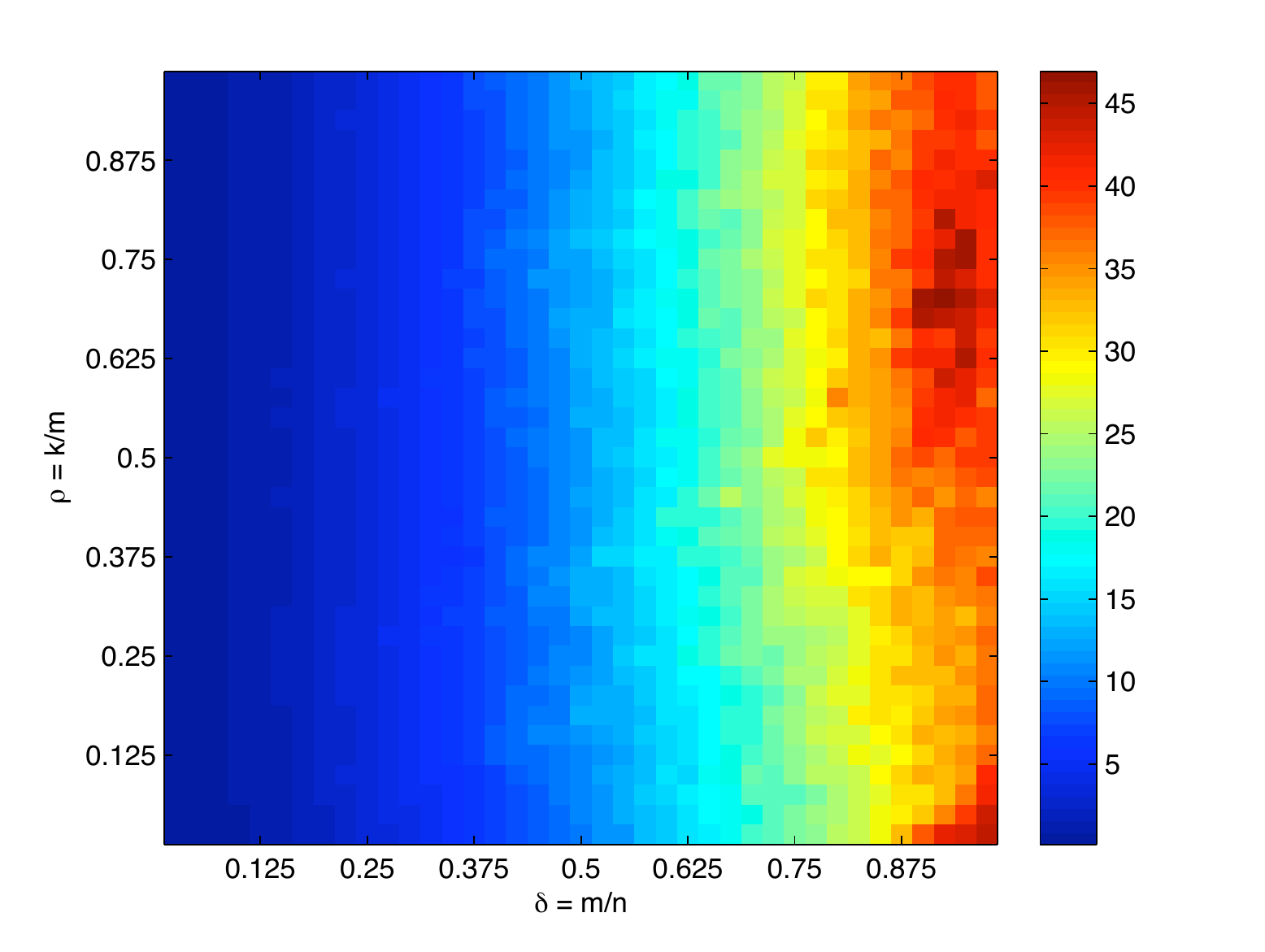}\end{center}
\caption{Duration of $\ell_1$ minimization (left) and solving \eqref{MinProg:DualCert} with injectivity check (right) in seconds.}
\label{Fig:Time_Comp_2}
\end{figure}

In dependence of $\delta$ and $\rho$, Figure \ref{Fig:Time_Comp} shows the averaged duration of solving \eqref{MinProg:DualCert} and calculating the rank of the submatrix divided by the averaged duration of the $\ell_1$ minimization.
One may observe that all quotients are less than one which means that in all cases solving \eqref{MinProg:DualCert} and checking the injectivity of the submatrix is faster than solving Basis Pursuit as a linear program.
Figure~\ref{Fig:Time_Comp_2} illustrates that the duration of both methods do increase with an increasing $\delta$, but while solving \eqref{MinProg:DualCert} seems to depend only on $\delta$, the $\ell_1$ minimization depends on $\delta$ and also on $\rho$.
Moreover, the contours of $\rho_F$ from Theorem \ref{Th:PropNeighbWeak} can be seen in the duration of time at the $\ell_1$ minimization as well as in the Figure \ref{Fig:Time_Comp_PT}: one may say that, on average, solving Basis Pursuit at $\rho=\rho_F(\delta)$ takes more time than solving it at any different $\rho$ in the neigborhood of $\rho_F(\delta)$.
Additionally, for small $\delta$ only small differences up to a quotient of $4/5$ appear in the comparision of the time duration.
In total, the use of checking~\eqref{MinProg:DualCert} instead of doing $\ell_1$ minimization reduces the computational time by a factor of $0.29$
(which amounts to a total save of 16 hours of computational time in our experiments).
% Indeed, most experiments  use a striktly higher number $M$ of repetitions, e.g. $M=200$ in \cite{DoTa09}, than we did.

Furthermore, one may observe that the quotients decrease between $\delta = 0.225$ and $\delta=0.25$.
We can not give reasonable causes for this phenomenon but remark that this process stems from the duration of the $\ell_1$ minimization program, cf. Figure~\ref{Fig:Time_Comp_2}.
%Additionally, Figure \ref{Fig:Time_Comp_2} shows another interesting phenomenon: the contours of $\rho_F$ from Theorem \ref{Th:PropNeighbWeak} can be seen in the duration of time at the $\ell_1$ minimization but not in the one for checking~\eqref{MinProg:DualCert}: On average, one may say that solving Basis Pursuit at $\rho=\rho_F(\delta)$ takes more time than solving it at any different $\rho$ in the neigborhood of $\rho_F(\delta)$.
% None of these surprising phenomenons can be observed in at the duaration of solving \eqref{MinProg:DualCert} with injectivity check, cf. Figure~\ref{Fig:Time_Comp_2} (right).

%This is highlighted by the fact that for fixed $\delta$ the value $\rho_F(\delta)$ is mostly larger than the rest of the values $\rho$ in this column.
%That means that for the pair $(\delta,\rho_F(\delta))$ our approach is once again remarkably faster than $\ell_1$ minimization.

\subsection{Number of Recoverable Supports for Certain Types of Matrices}\label{Se:NRSCM}
In this subsection, we compare computational experiments on the number of Recoverable Supports of several types of matrices with results from Section 3 whereas we restrict our experiments to minimally redundant matrices.
The computational experiments were done by Monte Carlo experiments described above.
Since in the previous sections only Gaussian matrices as well as Mercedes-Benz frames were considered, we will use these types as test problems.
Note that in any repetition of the Monte Carlo procedure, a new Gaussia matrix is drawn; the calculated value approximates the expected number of Recoverable Suppports.

Similar to Section \ref{Se:Com_l1_l_inf}, the experiments were done by checking \eqref{Conditions:l1} through checking whether the corresponding submatrix is injective and solving \eqref{MinProg:DualCert} afterwards.
If the optimal value is strictly less than $1-10^{-12}$, we record the chosen pair $(I,s)$ as a Recoverable Support.
We did the experiments with $n = 15, 34, 155$ and $n= 555$ and all $|I|=k\le n-1$.
For each $k$ we did $M=1000$ repetitions.

\begin{figure}
\begin{center}\includegraphics[width=100mm]{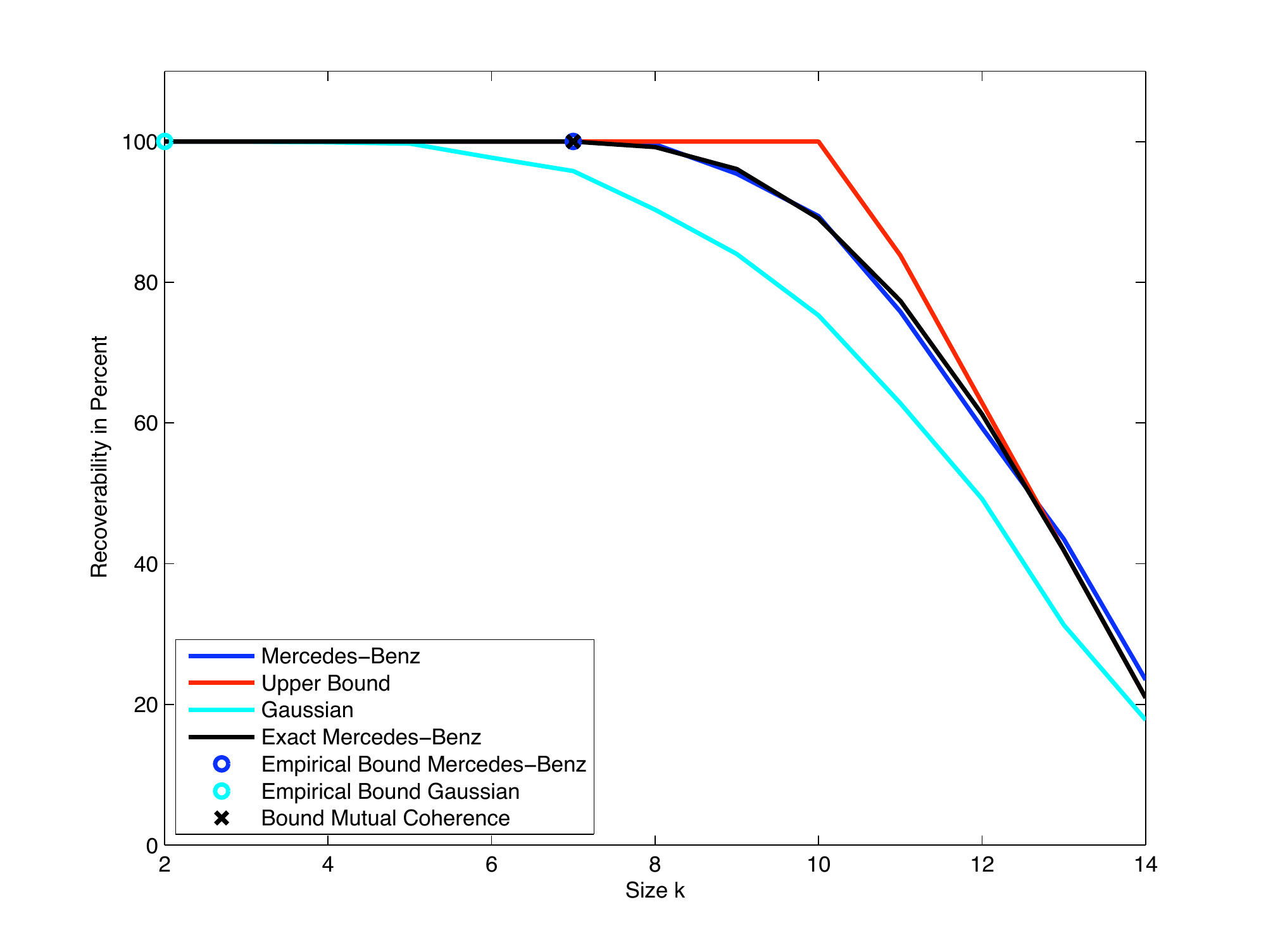}
\end{center}
\caption{Monte Carlo Sampling for $n=15$. The black curve respresents the actual number of Recoverable Supports of the Mercedes-Benz frame proportional to $2^k\binom{15}{k}$.}
\label{Fig:N=15}
\end{figure}

\begin{figure}
\begin{center}\includegraphics[width=60mm]{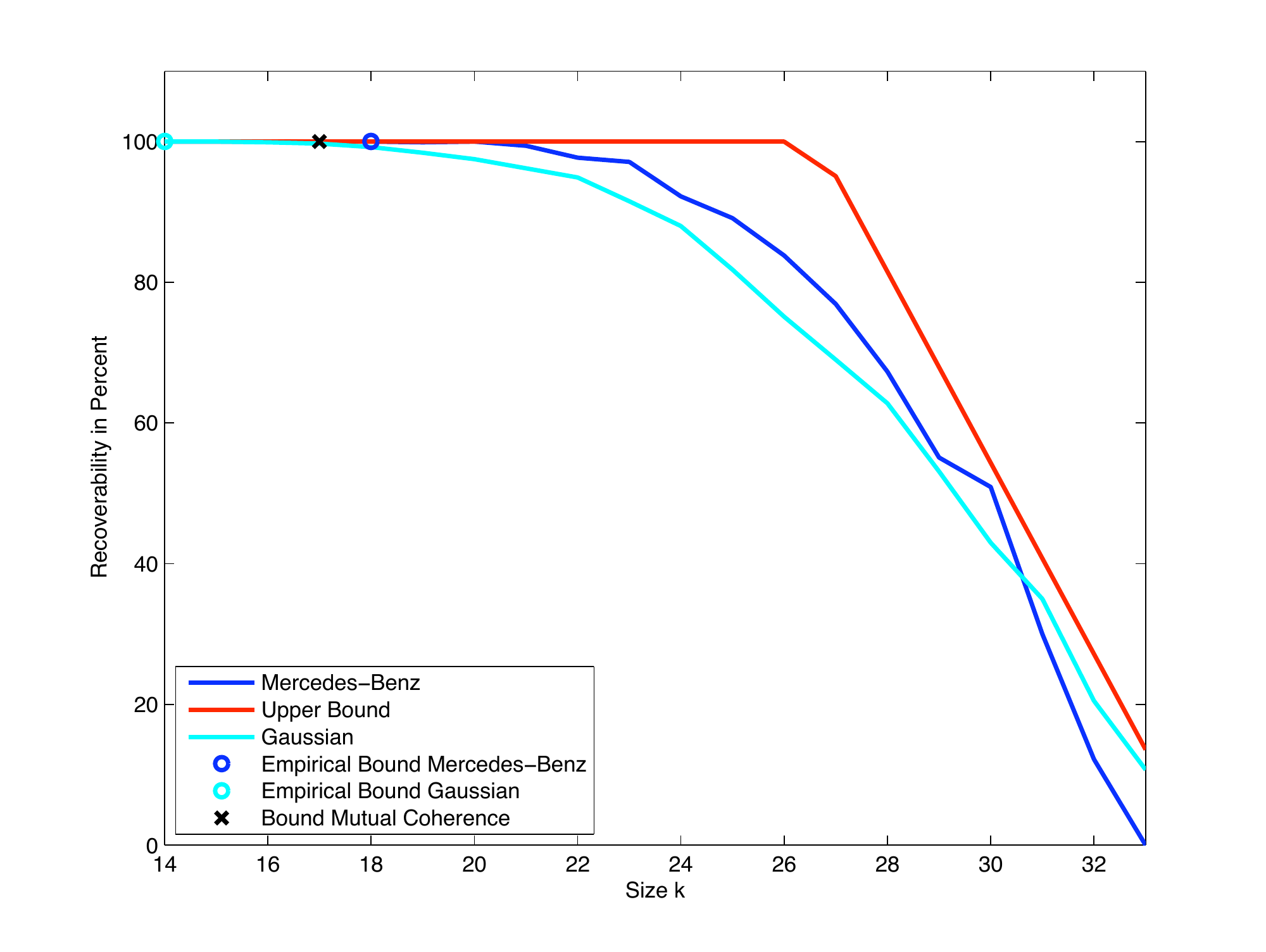}
\includegraphics[width=60mm]{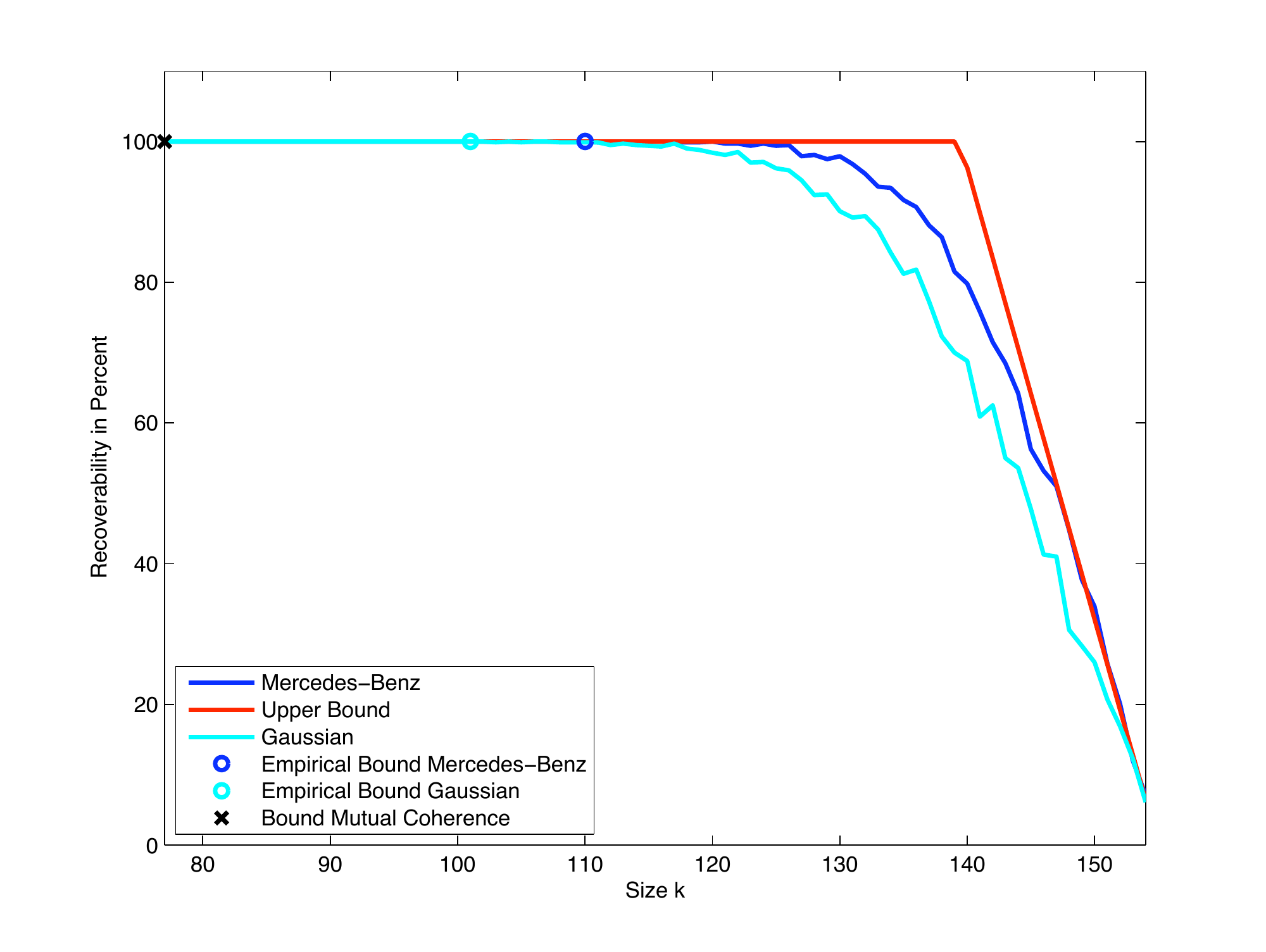}\end{center}
\caption{Monte Carlo Sampling for $n=34$ (left) and $n=155$ (right).}
\label{Fig:N=155}
\end{figure}

\begin{figure}
\begin{center}\includegraphics[width=60mm]{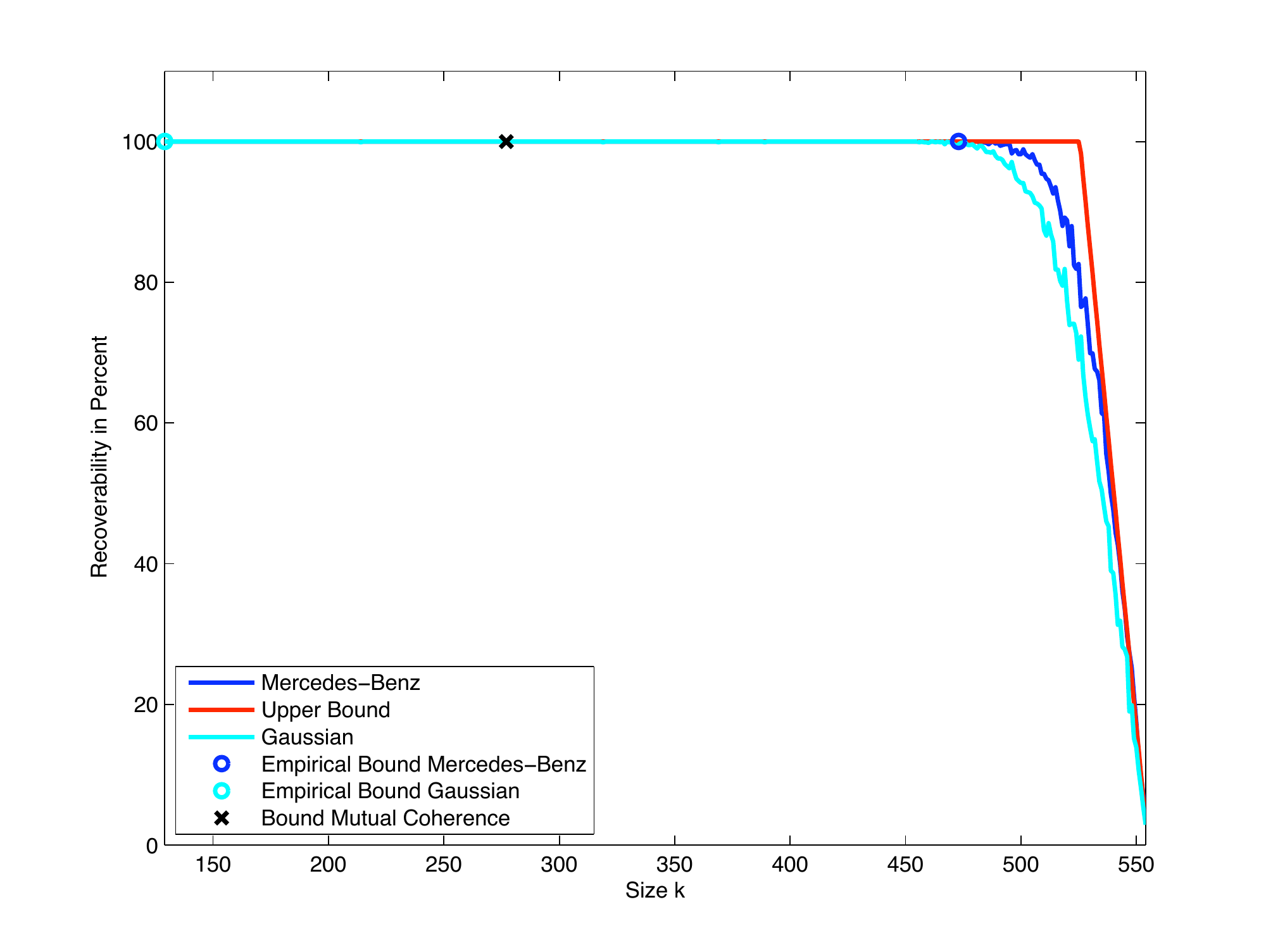}
\includegraphics[width=60mm]{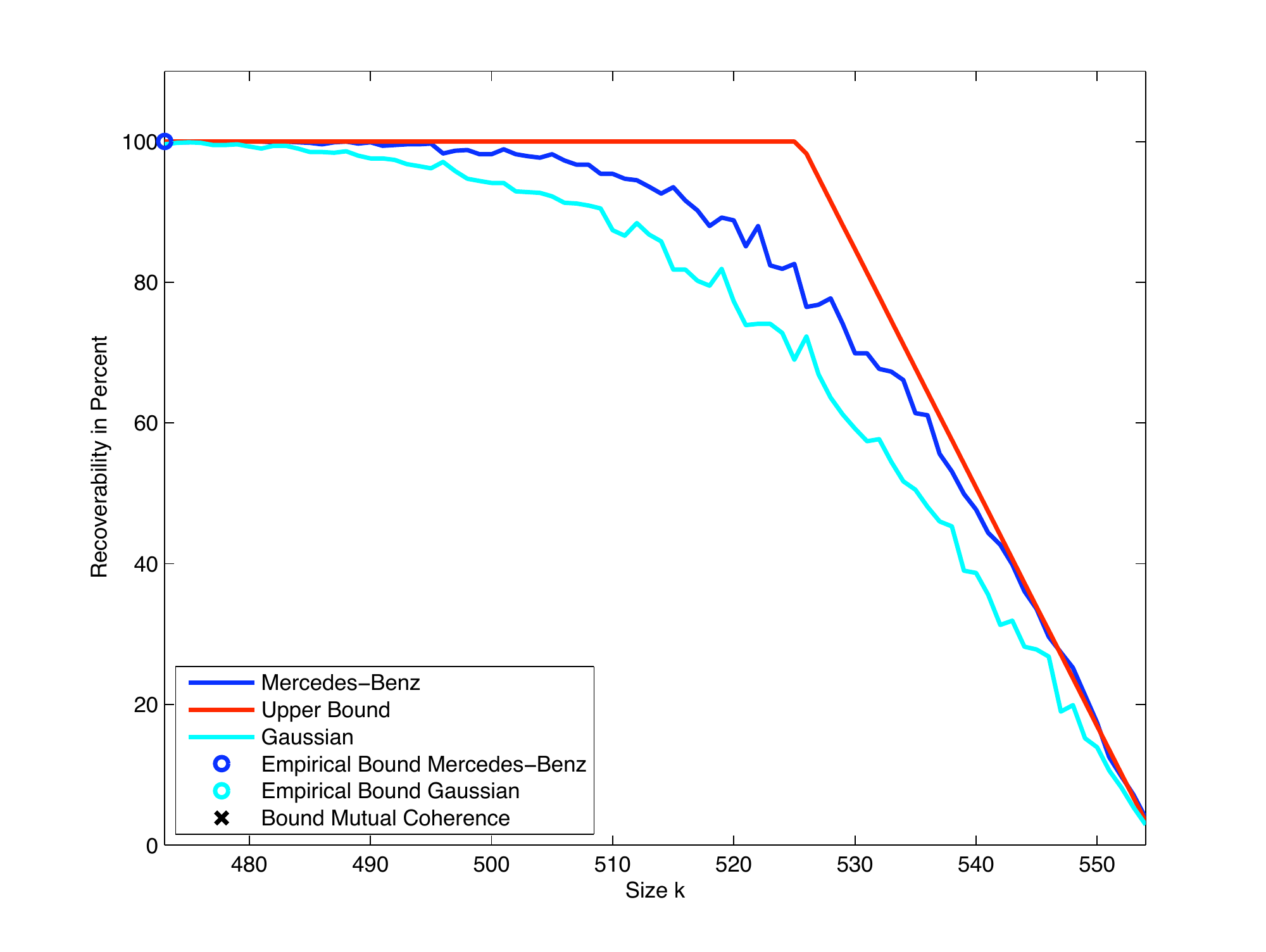}\end{center}
\caption{Monte Carlo Sampling for $n=555$. Left: segment from the ``Empirical Bound Gaussian'' to $k=554$. Right: segment from the ``Empirical Bound Mercedes-Benz'' to $k=554$, this graphics is a segment of the left graphics.}
\label{Fig:N=555}
\end{figure}

In Figures \ref{Fig:N=15}-\ref{Fig:N=555} all results are shown averaged.
The size $k$ of the desired Recoverable Support is given on the x-axis, on the y-axis the probability of recoverability is shown in percent.
These functions are empirical approximations of the mapping \eqref{Map:Recov}.
For comparison, the heuristic upper bound from \eqref{InEq:UpBoundBeta} in proportion to the total number $2^k\binom{n}{k}$ is also displayed.
Additionally, a circle for each type of matrix denotes the size $k$ when the recoverability at $k+1$ is less then one hundred percent (\textit{Empirical Bound}).
The empirical bounds are upper bounds for the smallest value $k$ where the actual recoverability \eqref{Map:Recov} at $k+1$ is less than one hundred percent, since there exists one pair $(I,s)$ which is not a Recoverable Support and the recoverability curve \eqref{Map:Recov} is monotonically nonincreasing by Proposition \ref{Pr:Recov}.
Note that in almost all cases (e.g. $n = 155, k = 111$ in Figure~\ref{Fig:N=155}) the empirical recoverability curves are not monotonically nonincreasing due their empirical nature.
The black cross denotes the last $k$ for which the recoverability guarantee for small sizes in Corollary \ref{Co:MuCo} holds.
All figures only show results from the smallest of all displayed bounds to $n-1$, since the tests deliver a recoverability of one hundred percent for the missing sizes.
Besides the empirical results for the Mercedes-Benz frame $A$, Figure \ref{Fig:N=15} shows the actual ratio $\Lambda(A,k)\left[2^k\binom{n}{k}\right]^{-1}$ in black with respect to $k$ for $n=15$.
For these results each of the $2^{|I|}\binom{15}{|I|}$ pairs $(I,s)$ with $|I|\le 14$ have been checked solving \eqref{MinProg:DualCert} if it was a Recoverable Support.


We emphasize that Mosek solved all problems with the status 'Optimal'.
In Figure \ref{Fig:N=15} one can see for the Mercedes-Benz frame that the results of the Monte Carlo sampling (blue) coincide with the actual values (black) up to an error of $10^{-1}$.
We tolerate this margin of error since improving the precision on one-tenth, we need to increase the number of samplings $M$ a hundredfold.
All results are bounded by the Upper Bound (red) except for the Mercedes-Benz frame in this case, which obviously is owed by the lack of accuracy.
Further the ``Bound Mutual Coherence'' coincides with the empirical bound for the Mercedes-Benz frame, which is not the case in the other cases.
Only in the case $n=155$ the ``Bound Mutual Coherence'' is the weakest bound, but as expected the distance to the ``Empirical Bound Mercedes-Benz''' increases with increasing $n$.
In all cases, Mercedes-Benz has the largest empirical bound.
At $n=155$, this values is $k=151$, while for $n=555$ it is $k=543$.
However, the distance between the 'Empirical Bound Mercedes-Benz' and the Upper Bound reaching one hundred percent increases with increasing $n$.
Additionally, Proposition \ref{Pr:SparsityK-1} holds for all suitable cases except an error of at most $10^{-2}$.
Hence, the results underlay the expectation that \eqref{InEq:UpBoundBeta} is a good bound for $k$ close to $n-1$.

Regarding Gaussian matrices, we observe that these matrices do not exceed the empirical recoverability curve of the Mercedes-Benz frame if $n$ is odd.
Contrary, it is expected that, at least with $k$ close to $n-1$, the recoverability curves of the Gaussian matrices exceed the curve of the Mercedes-Benz frame in case $n$ even; this behaviour may be observed in Figure~\ref{Fig:N=155}.
%Furthermore a unique differentiation of the Mean and the Single Gaussian matrices can not be noticed:
%For several $n$ the empirical bounds switch their succession, e.g. in $n=155$ the 'Empirical Bound Mean Gaussian' is the smallest of both bounds, while in $n=555$ it is larger than the 'Empirical Bound Single Gaussian'.
%Even a tendency to the distance to the 'Bound Mutual Coherence' is unapparent.

As also observed in the past similar experiments (e.g. \cite{TsDo06,BrDoEl09}), in all cases one can notice a rapid transition from one hundred to zero percent as $k$ increases.

%% file: section6_conclusion.tex
\section{Conclusion}

In this paper, we gave further insight in the apparently difficult question which vectors are recoverable by $\ell_1$ minimization for a given matrix $A$.
Through arranging recoverable vectors in equivalance classes (Recoverable Supports), dependent on $A$, it follows from Theorem~\ref{Th:poset} that the Recoverable Supports form a partial ordered set, which is completely known if its maximal elements, i.e. Maximal Recoverable Supports, are known.
Although Algortihm~\ref{Algo:ComputingRS} is able to compute such a Maximal Recoverable Support quite quickly, even for rather large matrices, we are still far away from any computational method which can result in an exhausting description of the set of Recoverable Supports (and such a method seems to be out of reach).

Moreover, we elaborated on a geometrical viewpoint on sparse recovery which is dual to the view through the projected cross polytope.
Exact values and new bounds on the number of Recoverable Supports were derived by connecting $\ell_1$ minimization to the dual approach via cross sections of the hypercube which has impact on probability whether a given vector can be reconstructed.

%%% Local Variables: 
%%% mode: latex
%%% TeX-master: "main"
%%% End: 